\documentclass[journal,twoside,web]{ieeecolor}
\usepackage{generic}
\usepackage{cite}
\usepackage{amsmath,amssymb,amsfonts}
\usepackage{graphicx}
\usepackage{algorithm,algpseudocode}
\usepackage{hyperref}
\hypersetup{hidelinks=true}
\usepackage{textcomp}
\usepackage{stmaryrd}
\usepackage{mathrsfs}
\usepackage{tikz}
\usetikzlibrary{shapes.geometric, arrows}
\graphicspath{{images/}}
\usepackage[colorinlistoftodos]{todonotes}
\usepackage{caption}
\usepackage{subcaption}
\usepackage{bm}
\usepackage[T1]{fontenc}
\usepackage{tcolorbox}
\usepackage{textcomp}
\usepackage{bigints}
\usepackage{cleveref}
\usepackage{siunitx}
\usepackage{textgreek}

\makeatletter
\algrenewcommand\alglinenumber[1]{\footnotesize \arabic{algorithm}.#1:}
\makeatother

\newcommand\restr[2]{{
		\left.\kern-\nulldelimiterspace 
		#1 
		\littletaller 
		\right|_{#2} 
}}

\newcommand{\littletaller}{\mathchoice{\vphantom{\big|}}{}{}{}}

\makeatletter
\newcommand{\algrule}[1][.2pt]{\par\vskip.5\baselineskip\hrule height #1\par\vskip.5\baselineskip}
\makeatother

\newtheorem{definition}{Definition}
\newtheorem{theorem}{Theorem}
\newtheorem{proposition}{Proposition}
\newtheorem{lemma}{Lemma}

\newtheorem{example}{Example}

\def\BibTeX{{\rm B\kern-.05em{\sc i\kern-.025em b}\kern-.08em
    T\kern-.1667em\lower.7ex\hbox{E}\kern-.125emX}}
\begin{document}
\title{Hierarchical Signal Temporal Logic Control of Continuous-Time Systems}
\author{Armin Pirastehzad, \IEEEmembership{Member, IEEE} and Bart Besselink \IEEEmembership{Senior Member, IEEE}
\thanks{This work was funded by the European Union (ERC, COCOS,101169999). Views and opinions expressed are however those of the authors only and do not necessarily reflect those of the European Union or the European Research Council. Neither the European Union nor the granting authority can be held responsible for them.}
\thanks{All authors are with the Bernoulli Institute for Mathematics, ComputerScience and Artificial Intelligence, University of Groningen, 9700 AK Groningen, The Netherlands (email: a.pirastehzad@rug.nl; b.besselink@rug.nl).}
}

\maketitle

\begin{abstract}
	For linear continuous-time systems, we solve the signal temporal logic (STL) control problem by translating it from continuous to discrete time such that a solution to the discrete-time problem guarantees a solution to the original continuous-time one. We split control synthesis into three hierarchical tasks. First, we discretize the continuous-time STL formula into a discrete-time one by sampling the temporal behavior it specifies. We accordingly discretize the continuous-time system by constructing a discrete-time model such that the input-state trajectories of the continuous-time system can be characterized as piecewise polynomial interpolations of the input-state trajectories of the discrete-time model. Second, we exploit this characterization to formulate mixed-integer linear constraints on the discrete-time trajectories that allow for planning a trajectory whose piecewise polynomial interpolation satisfies the continuous-time STL formula. In fact, system discretization and trajectory planning are performed in such a  manner that any feasible solution to these constraints gives a discrete-time trajectory that can be transformed into a continuous-time one guaranteed to satisfy the original STL formula. Third, we propose a step-by-step algorithm to construct the control input that realizes such continuous-time trajectory.  
\end{abstract}


\section{Introduction}\label{Sec_Introduction}
\IEEEPARstart{M}{odern} technologies, such as self-driving vehicles \cite{7823109}, power grids \cite{DOSSANTOSSERRA2025111625}, and healthcare robots \cite{Tavakoli2024}, often rely on timely execution of complex tasks under strict deadlines. These intricate tasks are typically specified in terms of signal temporal logic (STL) \cite{10.1007/978-3-540-30206-3_12}, which is a formal specification language that allows for description of time-dependent properties in system behavior. 

While control synthesis for STL formulas can be efficiently addressed for discrete-time systems, it remains a challenge for continuous-time ones. Inspired by the efficiency of control synthesis in discrete time, this paper aims to solve the continuous-time STL control problem by translating it into a discrete-time one such that a solution to the former can be obtained on the basis of a solution to the latter. This manifests in a \emph{hierarchical} scheme that splits the STL control task into 1) \emph{discretization}, where a discrete-time model of the system and a discretized version of the STL formula are derived, 2) \emph{discretization-aware planning}, where the discrete-time STL control problem (which corresponds to the derived model and STL formula) is solved, and 3) \emph{continuous-time execution}, through which a solution to the discrete-time STL control problem is transformed into a solution to the original continuous-time one.

\subsection{Related Works}
Existing solution schemes to the continuous-time STL control problem (\textit{e.g.,} those proposed in \cite{8404080,11036802,9992999,10388467,10144389,10421775,9147387,9161270,9483028,11505760,zaccherini2026controlmultiagentsystemsstl,zaccherini2026observerbasedcontrolmultiagentsystems}) usually have in common that they first translate the set of trajectories that adhere to the STL formula into (time-varying) controlled invariant sets and then seek the control input that establishes forward invariance. Specifically, these sets are specified as (super)level sets of  so-called (high-order/adaptive) time-varying control barrier functions \cite{7782377,8796030,GARG2024100945,9456981,9516971,XU2018195,9410332}, which allows for the characterization of forward invariance as a quadratic program whose solution gives the control input at each instant of time.

While these solution schemes are effective and mathematically elegant, they present a few major drawbacks. First, many of these schemes are limited in the sense that they can enforce only \emph{fragments} of STL formulas, \textit{i.e.,} they can only guarantee the satisfaction of a \emph{restricted} class of STL formulas. In particular, the schemes proposed in \cite{8404080,9161270,10421775} do not support \emph{nested} formulas (which are formed by composition of multiple temporal operators), whereas those proposed in \cite{11505760,zaccherini2026controlmultiagentsystemsstl,zaccherini2026observerbasedcontrolmultiagentsystems} support only specific types of nested formulas rather than arbitrary ones. Second, constraints on the control input encumber the implementation of these schemes, as they typically hinder the construction of a \emph{valid} control barrier function, which jeopardizes the feasibility of the quadratic program whose solution gives the appropriate control input (see, \textit{e.g.,} \cite{9683085}, for further details). Third, while in principle, this quadratic program allows for the explicit construction of the control input as a time-varying state feedback, in practice, such explicit construction is rarely possible. As a consequence, existing schemes often obtain the control input in a \emph{sampled-data} fashion, \textit{i.e.,} they obtain the control input at specific instants of time rather than at all times (see, \textit{e.g.,}  \cite[Section III.C]{8404080}, \cite[Section III.D]{10388467}). Implementation of such control, however, guarantees the STL formula only at these specific instants, while it remains inconclusive on its satisfaction in between them (see  \cite[Section I]{9147387} for a detailed discussion). Fourth, the efficacy of these schemes heavily relies on the choice of control barrier functions that encode the STL formula, \textit{i.e.,} some choices may render the problem infeasible even though feasible solutions exist (see, \textit{e.g.,} \cite[Section III.A]{8404080}). It is worth noting that some of these drawbacks are addressed by the solution schemes proposed in \cite{10388467,9147387}. Specifically, the approach presented in \cite{10388467} overcomes the first drawback by supporting arbitrary STL formulas, whereas the technique proposed in \cite{9147387} addresses the third drawback by providing guarantees for the inter-sample behavior. Nevertheless, both schemes suffer from the remaining drawbacks due to their reliance on (time-varying) control barrier functions.

While the continuous-time STL control problem remains a challenge, numerous efficient solution schemes have been proposed for its discrete-time counterpart, see, \textit{e.g.,} \cite{7039363,belta2019formal,8310901,9769752,9351643,9261610,9992326,YIN2024100940}. These schemes have in common that they exploit the \emph{discreteness} of solution trajectories to encode the STL formula as mixed-integer linear constraints on the control input sequence. This allows for utilization of efficient numerical algorithms and software (\textit{e.g.,} those listed in \cite{CLAUTIAUX2025707}) to solve a mixed-integer linear program (MILP) whose solution gives the control sequence that enforces the STL formula. Discrete-time solution schemes are broadly expressive as they support the full STL syntax. Additionally, these schemes are not susceptible to input constraints as they encode the STL formula in terms of an MILP, which naturally accommodates such constraints. More importantly, in contrast to the continuous-time STL control schemes that are conservative in the sense that they yield sufficient conditions for formula satisfaction (as they rely on the choice of control barrier functions that encode the formula), the discrete-time schemes are complete in the sense that they yield necessary and sufficient conditions for formula satisfaction. 

Motivated by these advantages, we develop a hierarchical scheme that solves the continuous-time STL control problem by translating it into a discrete-time one, whose solution can be efficiently obtained using the existing control schemes. We accomplish this by making the following contributions. 

\subsection{Contributions}
We propose a methodology for the translation of an STL control problem from continuous to discrete time such that a solution to the discrete-time problem \emph{guarantees} a solution to the original continuous-time one. As a matter of fact, we characterize a solution to the continuous-time STL control problem in terms of a solution to its discrete-time counterpart. Such characterization then leads to the derivation of a control input that enforces the continuous-time STL formula based on the control sequence that enforces the discrete-time STL formula. Our methodology proceeds as follows.

First, we discretize the continuous-time STL formula by sampling the temporal behavior it specifies and constructing a discrete-time STL formula that captures this sampled behavior. Then, using the same sampling time, we discretize the continuous-time system into a discrete-time model whose input-state behavior characterizes that of the original continuous-time system. Specifically, we conduct discretization according to the notion of \emph{system interpolation} \cite{Pirastzehzad2025} to construct this discrete-time model such that the input-state trajectories of the continuous-time system can be realized as piecewise polynomial interpolations of the input-state trajectories of the discrete model. We thus characterize the (input-state) behavior of the continuous-time system in terms of its discrete-time model. 

Second, we exploit this characterization to formulate mixed-integer linear constraints on the discrete-time input-state trajectories that allow for planning a trajectory whose piecewise polynomial interpolation, generated by the continuous-time system, satisfies the original STL formula. In other words, we perform system discretization and trajectory planning in such a manner that any feasible solution to these constraints gives a discrete-time trajectory that can be transformed into a continuous-time one guaranteed to satisfy the original STL formula. This leads to an MILP whose solution, after piecewise polynomial interpolation, gives the continuous-time control input that enforces the STL formula. As a consequence, in order to solve the continuous-time STL control problem, it suffices to solve an MILP that is completely in terms of the discrete-time model obtained from discretization by system interpolation. This is computationally advantageous as numerous efficient algorithms and software are available for solving such MILP (see, \textit{e.g.,} \cite{CLAUTIAUX2025707}, for further details). We accordingly propose a step-by-step algorithm for hierarchical continuous-time STL control. 

Our main contribution is therefore the mobilization of system interpolation to characterize a continuous-time STL control problem in terms of a discrete-time one such that a solution to the original continuous-time problem can be obtained by simply solving the discrete-time one, which can be addressed by numerous efficient control schemes. 

The hierarchical scheme proposed in this paper offers advantages over the existing solution schemes to the continuous-time STL control problem in the following aspects.

First, in contrast to solution schemes that address only fragments of STL formulas (\textit{i.e.,} restricted classes of formula structures), the proposed hierarchical scheme does not impose any assumptions on the formula structure. This is due to the fact that existing schemes typically rely on the characterization of STL formulas in terms of control barrier functions, which becomes prohibitively cumbersome for complex formula structures (\textit{e.g.,} those formed by nested temporal operators). The proposed scheme, however, relies on translation of the continuous-time STL control problem into a discrete-time one, which does not necessitate any a priori assumptions on the structure of STL formulas.  

Second, input constraints complicate the implementation of solution schemes that encode the STL formula in terms of control barrier functions, as these constraints typically hinder the construction of valid control barrier functions. In contrast, such complications do not arise in the proposed approach as it encodes the STL formula as an MILP, which naturally incorporates control input constraints. 

Third, solution schemes that obtain the control in a sampled-data fashion enforce the STL formula only at the sampled time instants rather than at all times, \textit{i.e.,} they do not necessarily guarantee the formula in between these instants. The proposed scheme, on the other hand, enforces the STL formula at all times, \textit{i.e.,} it provides guarantees for the formula satisfaction not only at sampled time instants but also over the time intervals between them. 

Fourth, solution schemes that encode STL formulas in terms of control barrier functions are inherently conservative as their efficacy critically relies on the choice of these functions, \textit{i.e.,} some choices may render the problem infeasible even though feasible solutions exist. Such conservationism is reduced in the proposed scheme as it encodes the STL formula as an MILP, which provides an equivalent representation of the STL formula at the sampled times. 

The rest of this paper is organized as follows. In Section~\ref{Sec_ProblemStatement}, we give a formal statement of the continuous-time STL control problem and, accordingly, propose the hierarchical solution scheme. In Section~\ref{Sec_Discretization}, we perform the discretization step of the hierarchical scheme, where we develop methodologies for discretization of the continuous-time system and its STL formula. We then proceed with the planning step of the hierarchical scheme in Section~\ref{Sec_Planning}, where we formulate the continuous-time STL control problem as an MILP that is completely in terms of the discrete-time model of the original system. In Section~\ref{Sec_CTExecution}, we perform continuous-time execution, which is the last step of the hierarchical scheme. Here, we propose a step-by-step algorithm for construction of the control input that enforces the continuous-time STL formula. We then illustrate our results in a numerical example in Section~\ref{Sec_Example}, and finally conclude the paper in Section~\ref{Sec_Conclusion}. 
 
 \subsection{Notation} 
 The mathematical notation adopted in this paper is largely standard. We briefly review the general conventions. 
 \subsubsection*{Sets}
 We respectively denote the set of real numbers and integers by $\mathbb{R}$ and $\mathbb{Z}$. We also denote the closed real interval $\{t\in\mathbb{R} \vert t_1\leq t \leq t_2\}$ and the closed integer interval $\{k\in\mathbb{Z} \vert k_1\leq k \leq k_2\}$ by $[t_1,t_2]$ and $\llbracket k_1,k_2 \rrbracket$, respectively. We define a convex polyhedron in $\mathbb{R}^n$ as a set of the form $\{x\in\mathbb{R}^n | Fx\leq g\}$, where $F\in\mathbb{R}^{m\times n}$ and $g\in\mathbb{R}^m$, and the inequality is understood componentwise. 
 \subsubsection*{Operators} We define the vectorization of a matrix $M =[M_1 \, M_2\,\cdots\,M_m]$, with columns $M_1,M_2,\cdots,M_m\in\mathbb{R}^n$, as $\operatorname{vec}(M) = (M_1^\top,M_2^\top,\cdots,M_m^\top)^\top$. We also utilize the symbol $\otimes$ to denote the Kronecker product. Indicating an identity matrix and a zero matrix of appropriate dimensions by $I$ and $0$, respectively, we define the Kronecker sum of square matrices $N\in\mathbb{R}^{n\times n}$ and $M\in\mathbb{R}^{m\times m}$ as $N\oplus M = N\otimes I + I\otimes M$. 
 \subsubsection*{Function Spaces and Operations}
We use $\mathbb{P}_N^n$ to denote the space of all $n$-dimensional vector-valued polynomials up to degree $N$.  Given a set $I$, we define the restriction of a function $f : I \rightarrow \mathbb{R}^n$ to the subset $I'\subset I$ as the function $\restr{f}{I'} : I' \rightarrow \mathbb{R}^n$ such that $\forall t \in I' : \restr{f}{I'} (t) = f(t)$. Then, for a set $\mathcal{F}$ of functions $f:I \rightarrow \mathbb{R}^n$, we define the restriction of $\mathcal{F}$ to $I'$ as 
 $\restr{\mathcal{F}}{I'} = \left\{\restr{f}{I'} \big\vert f\in\mathcal{F}\right\}$. 
\section{Problem Statement}\label{Sec_ProblemStatement}
We develop a \emph{hierarchical} control framework for enforcing STL formulas on continuous-time linear time-invariant systems. Towards that end, we first provide a brief overview of STL as a \emph{specification} language that describes the desirable temporal behavior of dynamical systems. We then formulate the continuous-time STL control problem of interest and, finally, propose the hierarchical approach that we adopt to solve it. 
\subsection{Signal Temporal Logic Specifications}\label{Sec_STL}
We specify the \emph{desirable} (temporal) behavior of continuous-time signals in terms of STL \cite{10.1007/978-3-540-30206-3_12}, which is a formal specification language built upon (atomic) predicate functions $\mu:\mathbb{R}^n \rightarrow \{\mathrm{True},\mathrm{False}\}$. Specifically, we consider linear predicates of the form
\begin{equation}\label{Predicate}
	\mu (s) := \left\{\begin{aligned}
		\mathrm{True}, \quad \text{if} \; \gamma ^\top s + \beta>0,\\
		\mathrm{False}, \quad \text{if} \; \gamma^\top s + \beta\leq 0,
	\end{aligned} \right.
\end{equation}
where $\gamma \in \mathbb{R}^n$ and $\beta \in \mathbb{R}$. An STL formula is recursively constructed according to the syntax
\begin{equation}\label{Syntax}
	\mathfrak{S} \!::= \!\mathrm{True}\big\vert\mu \big\vert \neg \mathscr{S} \big\vert \mathscr{S}^1 \wedge \mathscr{S}^2 \big\vert \operatorname{G}_{[t_1,t_2]}\mathscr{S} \big\vert \mathscr{S}^1 \operatorname{U}_{[t_1,t_2]} \mathscr{S}^2,
\end{equation}
where $\mu$ is a predicate given by \eqref{Predicate}, $\mathscr{S},\mathscr{S}^1,\mathscr{S}^2$ are STL subformulas, and $t_1,t_2\geq 0$ are such that $t_1\leq t_2 <\infty$. We note that $\neg$ and $\wedge$ are Boolean operators that are respectively read as \emph{negation} and \emph{conjunction}, whereas $\operatorname{G}$ and $\operatorname{U}$ are \emph{temporal} operators that are respectively read as \emph{always} (or \emph{globally}) and \emph{until}. 

Given a continuous-time signal $s: [0,T] \rightarrow \mathbb{R}^n$ (with sufficiently large $T>0$), for an STL formula $\mathfrak{S}$, we denote by $(s,t)\models \mathfrak{S}$ the validity of $\mathfrak{S}$ with respect to $s$ at time $t$. This validity is inductively determined in accordance with the following semantics:
	\begin{align}
		(s,t) & \models  \mu       \!\!\!\!&\Leftrightarrow & \ \gamma ^\top s(t) + \beta>0,\nonumber\\
		(s,t) & \models  \neg \mathscr{S} \!\!\!\!&\Leftrightarrow & \ \neg \big((s,t)\models \mathscr{S}\big),\label{Semantics}\\
		(s,t) &\models \mathscr{S}^1\wedge \mathscr{S}^2 \!\!\!\!& \Leftrightarrow & \ \big((s,t)\models \mathscr{S}^1\big) \wedge \big((s,t)\models \mathscr{S}^2\big),\nonumber\\
		(s,t) &\models \operatorname{G}_{[t_1,t_2]} \mathscr{S} \!\!\!\!& \Leftrightarrow & \ \forall \tau\! \in\!\![t+t_1,\!t+t_2]\!:\! (s,\tau)\models \mathscr{S},\nonumber\\
		(s,t)&\models \mathscr{S}^1\operatorname{U}_{[t_1,t_2]}\mathscr{S}^2 \!\!\!\!& \Leftrightarrow & \  \exists t'\!\in\!\![t+t_1,\!t+t_2]\!: \!\big(\!(s,t')\!\models \!\mathscr{S}^2\big)\nonumber\\
		&&& \ \wedge \big(\forall \tau \!\in \![t,t']\!:\!(s,\tau) \models \mathscr{S}^1\big).\nonumber
	\end{align}
	
Under these semantics, the signal $s$ is said to \emph{satisfy} the formula $\mathfrak{S}$ if $(s,0) \models \mathfrak{S}$. We denote this by $s\models\mathfrak{S}$ to emphasize that the formula $\mathfrak{S}$ is valid with respect to the \emph{entire} signal $s$ rather than to its evaluation from a specific time instant $t$ onward.
We also note that the above semantics further enable the characterization of other Boolean and temporal operators. Namely, the \emph{disjunction} operator, denoted by $\vee$, can be obtained in terms of negation and conjunction as $\mathscr{S}^1\vee\mathscr{S}^2 := \neg\big(\neg\mathscr{S}^1 \wedge \neg\mathscr{S}^2\big)$. On the other hand, the \emph{eventually} (or \emph{future}) operator, denoted by $\operatorname{F}$, is defined in terms of until as $\operatorname{F}_{[t_1,t_2]}\mathscr{S} := \mathrm{True} \operatorname{U}_{[t_1,t_2]} \mathscr{S}$ (further details on STL semantics can be found in, \textit{e.g.,} \cite[Chapter 3]{lindemann2025formal} and \cite[Section II]{7039363}).

We observe that any STL formula $\mathfrak{S}$ constructed according to the syntax \eqref{Syntax} is \emph{bounded-time} as it does not contain any temporal operators with unbounded time intervals. We define the \emph{horizon} of $\mathfrak{S}$, denoted by $\mathfrak{h}_{\mathfrak{S}}$, as the time length required to evaluate the satisfaction of $\mathfrak{S}$, and we note from the semantics \eqref{Semantics} that $\mathfrak{h}_{\mathfrak{S}} < \infty$. We highlight that the validation of $\mathfrak{S}$ with respect to a continuous-time signal $s:[0,T] \rightarrow \mathbb{R}^n$ (with $T>\mathfrak{h}_\mathfrak{S}$) depends only on the behavior of $s$ over the time interval $[0,\mathfrak{h}_{\mathfrak{S}}]$. We thus observe that 
\begin{equation}\label{bounded}
	s\models \mathfrak{S} \Longleftrightarrow \restr{s}{[0,\mathfrak{h}_{\mathfrak{S}}]}\models \mathfrak{S},
\end{equation}
see, \textit{e.g.,} \cite[Section 3]{belta2019formal} for further details. 

With a few minor adjustments, STL can be also used to specify the temporal behavior of discrete-time signals. By replacing real-valued (time) intervals with integer-valued ones in \eqref{Syntax} and \eqref{Semantics}, we respectively derive the \emph{discrete-time} STL syntax and semantics. We emphasize that the definitions and notation for  discrete-time STL formulas are entirely similar to those in the continuous-time case. In particular, the bounded-time property of a discrete-time STL formula and its horizon are defined in exactly the same fashion as in the continuous-time case, which implies that a discrete-time counterpart to the equivalence \eqref{bounded} also holds.  

\subsection{Continuous-time Signal Temporal Logic Control}\label{Sec_ControlProblem}
Consider the continuous-time linear system 
\begin{equation}\label{ContinuousSystem}
	\bm{\Sigma}_c: \begin{aligned}
		\dot{x}_c(t) &= A_cx_c(t) + B_cu_c(t),
	\end{aligned}
\end{equation}
with state $x_c\in\mathbb{R}^n$ and control input $u_c \in\mathbb{R}^m$. We denote by $x_c(t;x_{0},u_c)$ the state solution, at time $t$, of \eqref{ContinuousSystem} for initial condition $x_c(0) = x_{0}$ and input function $u_c$.

Given an initial condition $x_0\in\mathbb{R}^n$, we are interested in constructing a continuous-time control input $u_c$ to enforce an STL formula $\mathfrak{S}_c$. This, as a consequence of \eqref{bounded}, is equivalent to 
\begin{equation}\label{ControlProblem}
	\restr{x_c(\cdot;x_0,u_c)}{[0,\mathfrak{h}_{\mathfrak{S}_c}]} \models \mathfrak{S}_c,
\end{equation}
which indicates that it suffices to design the control input over the finite interval $[0,\mathfrak{h}_{\mathfrak{S}_c}]$. The designed control input, however, must also satisfy constraints arising from physical and/or operational limitations imposed on the system. We model these constraints by restricting the control input to take values in a \emph{convex} polyhedral set $\mathcal{U}\subset\mathbb{R}^m$, \textit{i.e.,} $u_c(t) \in \mathcal{U}$ for all $t\geq 0$. Hence, we aim to solve the following continuous-time STL control problem. 

\textit{Continuous-time STL Control Problem:} Given a continuous-time control system \eqref{ContinuousSystem} and an STL formula $\mathfrak{S}_c$ constructed according to the syntax \eqref{Syntax}, for an initial condition $x_0\in\mathbb{R}^n$ and a convex polyhedral set $\mathcal{U}\subset\mathbb{R}^m$, find a continuous-time control input $u_c : [0,\mathfrak{h}_{\mathfrak{S}_c}]\rightarrow\mathcal{U}$ such that \eqref{ControlProblem} holds. 

Existing solution schemes to the continuous-time STL control problem (\textit{e.g.,} those proposed in \cite{8404080,11036802,9992999,10388467,10144389,10421775,9147387,9161270,9483028,11505760,zaccherini2026controlmultiagentsystemsstl,zaccherini2026observerbasedcontrolmultiagentsystems}) commonly encode the STL formula in terms of (high-order/adaptive) time-varying control barrier functions \cite{7782377,8796030,GARG2024100945,9456981,9516971,XU2018195,9410332}, which leads to the formulation of a quadratic program whose solution gives the control input that enforces the formula. These schemes, however, suffer from at least one the following issues: 
\begin{itemize}
	\item Limited expressiveness, as many of these schemes support only a subset of the STL syntax;
	\item Susceptibility to input constraints, which results from the challenges associated with the construction of valid control barrier functions in the presence of such constraints;
	\item Conservatism, since the effectiveness of these schemes heavily relies on the choice of control barrier functions;
	\item Lack of guarantee for inter-sample behavior, since these schemes practically obtain the control input in a sampled-data fashion and, therefore, guarantee the STL formula only at the sampled times, while they provide no guarantees over intervals between the sampled times. 
\end{itemize}

Solution schemes to the \emph{discrete-time} STL control problem, on the other hand, have in common that they encode the STL formula as an MILP. These schemes are broadly expressive as they support the full STL syntax. Furthermore, they can easily handle control input constraints as the MILP formulation naturally accommodates such constraints. More importantly, these schemes are complete as this MILP is an equivalent representation of the STL formula, \textit{i.e.,} the feasibility of this MILP is necessary and sufficient for the formula satisfaction. 

Motivated by these observations, we develop a hierarchical scheme that translates the continuous-time STL control problem into a discrete-time one such that a solution to the original continuous-time problem can be characterized in terms of a solution to its discrete-time counterpart. In fact, we design this translation in such a way that 1) discrete-time STL control techniques can be utilized to enforce STL formulas constructed from the \emph{full} syntax \eqref{Syntax} rather than a restricted subset, 2) satisfaction of the STL formula is guaranteed not only at sampled times but also over intervals between them, and 3) control input constraints can be handled efficiently.  

\subsection{Hierarchical Continuous-time STL Control}\label{Sec_HierarchicalScheme}
Given an STL formula $\mathfrak{S}_c$ constructed according to the syntax \eqref{Syntax}, we will \emph{hierarchically} synthesize a continuous-time control input $u_c:[0,\mathfrak{h}_{\mathfrak{S}_c}]\rightarrow\mathcal{U}$ that establishes \eqref{ControlProblem} according to the following consecutive steps.

\begin{figure}
	\centering
	\begin{tikzpicture}[scale=0.5]
		\node  (ContSys) at (0,0) [rectangle, 
		minimum width=0.3cm, 
		minimum height=0.5cm,
		text centered, 
		draw=black, 
		line width = 0.5pt,
		fill=gray! 10!, align=center] {\fontsize{10}{30} \selectfont $\bm{\Sigma}_c$};
		
		\node (uc) at (-3,0) {\fontsize{10}{30} \selectfont $u_c$};
		\node (xc) at (3,0) {\fontsize{10}{30} \selectfont $x_c$};
		\node (mod) at (4.5,0) {\fontsize{10}{30} \selectfont $\models$};
		\node (Phic) at (6,0) {\fontsize{10}{30} \selectfont $\mathfrak{S}_c$};
		\draw[line width = 0.25pt,->,>=stealth] (uc)--(ContSys);
		\draw[line width = 0.25pt,->,>=stealth] (ContSys)--(xc);
		\node  (DiscSys) at (0,-5) [rectangle, 
		minimum width=0.3cm, 
		minimum height=0.5cm,
		text centered, 
		draw=black, 
		line width = 0.5pt,
		fill=gray! 10!, align=center] {\fontsize{10}{30} \selectfont $\bm{\Sigma}_d$};
		\draw[line width = 0.25pt,->,>=stealth, dashed] (ContSys)--(DiscSys);
		\node  [rotate = 90] (exec) at (-0.3,-2.5) {\fontsize{6}{30} \selectfont Discretization};
		
		\node (ud) at (-3,-5) {\fontsize{10}{30} \selectfont $u_d$};
		\node (xd) at (3,-5) {\fontsize{10}{30} \selectfont $x_d$};
		\node (modd) at (4.5,-5) {\fontsize{10}{30} \selectfont $\models$};
		\node (Phid) at (6,-5) {\fontsize{10}{30} \selectfont $\mathfrak{S}_d$};
		\draw[line width = 0.25pt,->,>=stealth] (ud)--(DiscSys);
		\draw[line width = 0.25pt,->,>=stealth] (DiscSys)--(xd);
		
		\draw[line width = 0.25pt,->,>=stealth, dashed] (Phic)--(Phid);
		\node  [rotate = 90] (exec) at (5.7,-2.5) {\fontsize{6}{30} \selectfont Discretization};
		
		\draw[line width = 0.25pt,->,>=stealth, dashed] (ud)--(uc);
		\node  [rotate = 90,text width=2cm] (exec) at (-3.5,-2.1) {\fontsize{6}{30} \selectfont Continuous-time \\[-7pt] \hspace*{12pt}Execution};
		
		\draw[line width = 0.7pt,-,dotted] (-3.5,-5.75) --(7,-5.75) -- (7,-4.1) -- (-3.5,-4.1) -- (-3.5,-5.75);
		\node (planing) at (1.75,-6.25) {\fontsize{6}{30} \selectfont Planing};
	\end{tikzpicture}
	\caption{We split the continuous-time STL control into 1) discretization, where the STL formula $\mathfrak{S}_c$ and the system $\bm{\Sigma}_c$ are respectively discretized into an STL formula $\mathfrak{S}_d$ and a discrete-time model $\bm{\Sigma}_d$, 2) planning, through which a control sequence $u_d$ is planned to enforce the STL formula $\mathfrak{S}_d$ on the discrete-time model $\bm{\Sigma}_d$, and 3) continuous-time execution, which interpolates the planned control sequence $u_d$ into the continuous-time input $u_c$ such that \eqref{ControlProblem} holds.}\label{Hierarchical} 
\end{figure}
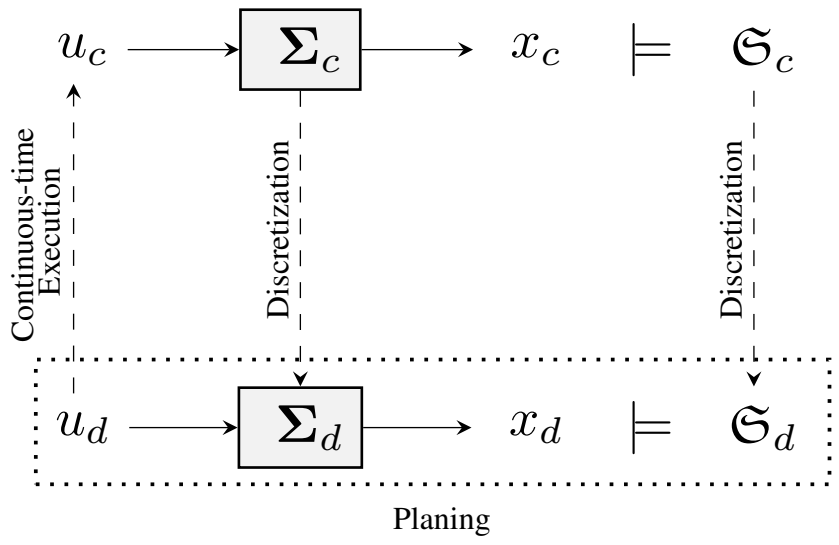

\begin{enumerate}
	\item \textit{Discretization}: For some suitable sampling time, we \emph{discretize} the continuous-time STL formula $\mathfrak{S}_c$ to obtain a (bounded-time) \emph{discrete-time} STL formula $\mathfrak{S}_d$ (as illustrated in the right sides of Figure~\ref{Hierarchical}). With respect to the same sampling time, we then \emph{discretize} the continuous-time system $\bm{\Sigma}_c$ into a discrete-time model $\bm{\Sigma}_d$ (as demonstrated in the middle part of Figure~\ref{Hierarchical}) in such a way that the input-state trajectories of $\bm{\Sigma}_c$ can be realized as piece-wise polynomial interpolations of the input-output trajectories of $\bm{\Sigma}_d$. This yields a characterization of the input-state behavior of $\bm{\Sigma}_c$ in terms of $\bm{\Sigma}_d$. Such characterization later enables the translation of specifications on $\bm{\Sigma}_c$ into those on $\bm{\Sigma}_d$. 
	
	We construct the discrete-time model $\bm{\Sigma}_d$ as a system of the form
	\begin{equation}\label{DiscreteSystem}
		\bm{\Sigma}_d: \begin{aligned}
			x_d(i+1)&=A_dx_d(i) + B_du_d(i),
		\end{aligned}
	\end{equation}
	with state $x_d\in\mathbb{R}^n$ and $u_d\in\mathbb{R}^m$. We use notation similar to that of \eqref{ContinuousSystem} to denote the state sequence of \eqref{DiscreteSystem}. 
	\item \textit{Discretization-aware Planning}: We \emph{plan} a discrete-time control sequence $u_d:\llbracket0,\mathfrak{h}_{\mathfrak{S}_d}\rrbracket\rightarrow\mathcal{U}$ such that 
	\begin{equation}\label{DiscreteControlProblem}
		\restr{x_d(\cdot; x_0, u_d)}{\llbracket0,\mathfrak{h}_{\mathfrak{S}_d}\rrbracket} \models \mathfrak{S}_d.
	\end{equation}
	This is a \emph{discrete-time} STL control problem (depicted in the lower part of Figure~\ref{Hierarchical}), which can be equivalently formulated as an MILP in terms of $u_d$. This MILP can be efficiently solved by software packages such as MATLAB's optimization toolbox and YALMIP \cite{Lofberg2004}. 
	
	\item \textit{Continuous-time Execution}: We exploit the fact that input-state trajectories of $\bm{\Sigma}_c$ can be realized as piece-wise polynomial interpolations of input-state trajectories of $\bm{\Sigma}_d$ to construct $u_c$ such that \eqref{ControlProblem} holds (see the left side of Figure~\ref{Hierarchical}). 
\end{enumerate}

Our solution scheme is therefore to \emph{translate} the continuous-time STL control problem into a discrete-time one whose solutions characterize those of the original (continuous-time) problem. In fact, we translate specifications on the continuous-time system into those on its discrete-time model in such a way that the satisfaction of these discrete-time specifications \emph{guarantees} the satisfaction of the original continuous-time ones. 

\section{Discretization}\label{Sec_Discretization}
We propose a methodology for discretization of the STL formula $\mathfrak{S}_c$ and the continuous-time system $\bm{\Sigma}_c$ that later allows for translation of the continuous-time STL control problem into a discrete-time one. 

\subsection{Formula Discretization}\label{Sec_SpecificationDiscretization}
Consider a continuous-time STL formula $\mathfrak{S}_c$. We discretize $\mathfrak{S}_c$ by sampling the temporal behavior it specifies to construct a discrete-time STL formula $\mathfrak{S}_d$ that preserves this temporal behavior at the sampled times.

We select a sampling time $\tau>0$ such that every time instant specified by $\mathfrak{S}_c$ (\textit{i.e.,} every endpoint of the time intervals in $\mathfrak{S}_c$) is an integer multiple of $\tau$. For this sampling time, we draw inspiration from the syntax \eqref{Syntax} to propose a \emph{recursive} procedure for construction of the discrete-time formula $\mathfrak{S}_d$. Given a (continuous-time) $\mathfrak{S}_c$, we construct $\mathfrak{S}_d$ by recursively \emph{transforming} subformulas of $\mathfrak{S}_c$ into their discrete-time counterparts according to the following procedure:
\begin{equation}\label{SpecificationDiscretizationProcedure}
 \begin{aligned}
 	\mu_c = \mu &\longrightarrow \mu_d = \mu,\\ 
 	\neg \mathscr{S}_c &\longrightarrow \neg \mathscr{S}_d,\\
 	\mathscr{S}_c^1 \wedge \mathscr{S}_c^2 &\longrightarrow \mathscr{S}_d^1 \wedge \mathscr{S}_d^2,\\
 	\mathscr{S}_c^1 \vee \mathscr{S}_c^2 &\longrightarrow \mathscr{S}_d^1 \vee \mathscr{S}_d^2,\\
 	\operatorname{G}_{[t_1,t_2]}\mathscr{S}_c &\longrightarrow \operatorname{G}_{\llbracket \frac{t_1}{\tau},\frac{t_2}{\tau}\rrbracket}\mathscr{S}_d,\\
 	\operatorname{F}_{[t_1,t_2]}\mathscr{S}_c &\longrightarrow \operatorname{F}_{\llbracket \frac{t_1}{\tau},\frac{t_2}{\tau}\rrbracket}\mathscr{S}_d,\\
 	\mathscr{S}_c^1 \operatorname{U}_{[t_1,t_2]}\mathscr{S}_c^2 &\longrightarrow \mathscr{S}_d^1 \operatorname{U}_{\llbracket \frac{t_1}{\tau},\frac{t_2}{\tau}\rrbracket}\mathscr{S}_d^2,
 \end{aligned}  
\end{equation}
where $\mu_c$ is a predicate in $\mathfrak{S}_c$ whose discrete-time counterpart in $\mathfrak{S}_d$ is $\mu_d$, whereas $\mathscr{S}_c$, $\mathscr{S}_c^1$, and $\mathscr{S}_c^2$ are subformulas in $\mathfrak{S}_c$ whose discrete-time counterparts in $\mathfrak{S}_d$ are respectively $\mathscr{S}_d$, $\mathscr{S}_d^1$, and $\mathscr{S}_d^2$. We observe from the recursive procedure \eqref{SpecificationDiscretizationProcedure} that 
\begin{equation}\label{Horizons}
	\mathfrak{h}_{\mathfrak{S}_c} = \mathfrak{h}_{\mathfrak{S}_d} \tau. 
\end{equation}
We demonstrate this procedure in the following example. 
\begin{example}\label{Example1}
	Consider the continuous-time STL formulas
	\begin{equation*}
		\begin{aligned}
			\mathfrak{S}_c^1 &= \operatorname{F}_{[0.5, 0.7]} \mu_1(s) \wedge \operatorname{G}_{[1,2]}\mu_2(s),\\
			\mathfrak{S}_c^2 &= \operatorname{F}_{[0.5, 1]}\operatorname{G}_{[3,3.5]} \mu(s),
		\end{aligned}
	\end{equation*}
	where $\mu, \mu_1,\mu_2$ are predicates of the form \eqref{Predicate}. As for $\mathfrak{S}_c^1$, we let $\tau = 0.1$ and apply the procedure \eqref{SpecificationDiscretizationProcedure} to discretize $\mathfrak{S}_c^1$ into the discrete-time STL formula $\mathfrak{S}_d^1$. Bearing in mind that predicates and Boolean operators remain intact and that only the temporal ones are affected by discretization, we obtain $\mathfrak{S}_d^1 = \operatorname{F}_{\llbracket5, 7\rrbracket} \mu_1(s)\wedge\operatorname{G}_{\llbracket10,20\rrbracket}\mu_2(s)$. Similarly, considering $\mathfrak{S}_c^2$, we let $\tau = 0.5$ and obtain $\mathfrak{S}_d^2 = \operatorname{F}_{\llbracket1,2\rrbracket}\operatorname{G}_{\llbracket6,7\rrbracket}\mu(s)$. 
	\end{example}
\subsection{System Discretization}\label{Sec_SystemDiscretization}
We discretize $\bm{\Sigma}_c$ according to the notion of system interpolation \cite{Pirastzehzad2025}, which determines whether the input-state trajectories of a continuous-time system can be realized as piecewise polynomial interpolations of the input-state trajectories of a discrete-time system. Towards that end, we first briefly review basic properties of shifted \emph{Legendre} polynomials, which we then exploit to define and characterize system interpolation. We accordingly use this characterization to discretize the continuous-time system $\bm{\Sigma}_c$ into a discrete-time model $\bm{\Sigma}_d$.
\subsubsection{Legendre Polynomials}\label{Sec_Legendre}
Given any integer $k\geq 0$, let $\mathfrak{L}_k^\tau$ denote the \emph{shifted} Legendre polynomial of degree $k$, over the interval $[0,\tau]$, defined as
\begin{equation}\label{ShiftedLegendre}
	\mathfrak{L}_k^\tau(t) = \frac{1}{2^k} \sum_{i = 0}^{\lfloor \frac{k}{2}\rfloor} (-1)^i \begin{pmatrix}
		k\\i
	\end{pmatrix} \begin{pmatrix}
		2k-2i\\k
	\end{pmatrix} \left(\frac{2t}{\tau} -1\right)^{k-2i},
\end{equation}
where $\lfloor \frac{k}{2}\rfloor$ is the greatest integer less than or equal to $\frac{k}{2}$.

We utilize the shifted Legendre polynomials to obtain a set of nodes such that any polynomial function, over the interval $[0,\tau]$, can be represented completely in terms of its values at these nodes. Considering polynomials up to degree $N\geq 1$, we choose the nodes $0 = t_{N,0}<t_{N,1}<\ldots<t_{N,N} < \tau$ such that 
\begin{equation}\label{InterpolationNodes}
	\forall i\in\{0,1,\cdots,N\}: \quad \mathfrak{L}_N^\tau(t_{N,i}) + \mathfrak{L}_{N+1}^\tau(t_{N,i}) = 0,
\end{equation}
\textit{i.e.,} we choose $t_{N,0}$, $t_{N,1}$, $\ldots$, $t_{N,N}$ as the zeros of $\mathfrak{L}_N^\tau + \mathfrak{L}_{N+1}^\tau$. We emphasize that such nodes always exists (further details can be found in \cite[Section 2.2.3]{Canuto2006} and \cite[Section 2]{wang2012legendre}). Based on these nodes, we now construct the polynomial function
\begin{equation}\label{PhiFunction}
	\phi_i(t) = \prod_{j=0,j\neq i}^{N}\frac{t-t_{N,j}}{t_{N,i}-t_{N,j}},
\end{equation}
for $i=0,1,\ldots,N$ (see, \textit{e.g.,} \cite[Section 3.1]{faires2012numerical}). We then note that 
\begin{equation}\label{Dirac}
	\phi_i(t_{N,j}) = \left\{\begin{aligned}
		&0&j&\neq i,\\
		&1&j&=i,
	\end{aligned} \right.
\end{equation}
which suggests that the polynomial functions $\phi_0, \phi_1,\ldots,\phi_N$ are linearly independent and, therefore, form a basis for $\restr{\mathbb{P}_N^{n}}{[0,\tau]}$. In fact, it can be shown (see, \textit{e.g.,} \cite[Theorem 3.2]{faires2012numerical}) that, for any polynomial $P$ of degree $N$,  
\begin{equation}\label{SpectralRepresentation}
	\forall t\in[0,\tau]: \quad P(t) = \sum_{i=0}^{N} \phi_i(t)P(t_{N,i}),
\end{equation}
\textit{i.e.,} the functions $\phi_i$ in \eqref{SpectralRepresentation} allow for expressing a polynomial $P$ completely in terms of its values at the nodes. We emphasize that the functions $\phi_i$ are \emph{independent} of $P$. 

\subsubsection{System Interpolation}\label{Sec_SystemInterpolation}
System interpolation is a comparative framework that inspects the ability of a continuous-time system $\bm{\Sigma}_c$ in reproducing piecewise polynomial interpolations of the input-state trajectories of a discrete-time system $\bm{\Sigma}_d$. The formal definition of this notion, however, relies on the following description for interpolation of discrete-time signals into continuous-time piecewise polynomial ones.
\begin{definition}\label{DefSignalInterpolation}
	For integers $\ell,p\geq 1$, let a discrete-time signal $s_d: \llbracket0,\ell \rrbracket \rightarrow \mathbb{R}^p$ be given. For a sampling time $\tau>0$ and an integer $N\geq 1$, a continuous-time signal $s_c: [0,\ell \tau]\rightarrow\mathbb{R}^p$ is said to be an $N$-th order interpolation of $s_d$ with respect to $\tau$ if there exist polynomial functions $P_0,P_1,\cdots,P_{\ell-1}\in\restr{\mathbb{P}_N^{p}}{[0,\tau]}$ such that for all $i \in \llbracket0,\ell-1\rrbracket$, we have that 
	\begin{subequations}\label{InterpolationCondition}
		\begin{equation}\label{InstanceEquivalence}
			\begin{aligned}
				s_c(i\tau) &= s_d(i),\\
				s_c\big((i+1)\tau\big) &= s_d(i+1)
			\end{aligned}
		\end{equation}
		\text{and}
		\begin{equation}\label{PolynomialRestriction}
			\forall t\in[0,\tau]: \quad s_c(i\tau+t) = P_i(t).
		\end{equation}
	\end{subequations}
\end{definition}
\vspace*{2mm}

We thus emphasize that an $N$-th order interpolation of a discrete-time signal $s_d: \llbracket0,\ell \rrbracket \rightarrow \mathbb{R}^p$ (with respect to sampling time $\tau$) is any continuous-time signal $s_c: [0,\ell \tau]\rightarrow\mathbb{R}^p$ whose values at any time $i\tau$ equals the value of $s_d$ at instance $i$, while its restriction to the interval $[i\tau,(i+1)\tau]$ can be described by a polynomial in $\restr{\mathbb{P}_N^{p}}{[0,\tau]}$. Considering that such continuous-time signal is not necessarily unique (for $N\geq 2$), we denote by $\mathbb{I}_{N}^\tau(s_d)$ the set of all $N$-th order interpolations of $s_d$ with respect to $\tau$. 

Based on this, we now give the formal definition of system interpolation, which is taken from \cite[Definition 2]{Pirastzehzad2025}.

\begin{definition}\label{DefInterpolatingSystem}
	For a sampling time $\tau>0$ and an integer $N\geq 1$, a continuous-time system $\bm{\Sigma}_c$ is said to be an $N$-th order interpolator of a discrete-time system $\bm{\Sigma}_d$ with respect to $\tau$ if for any integer $\ell\geq 1$, any initial condition $x_0\in\mathbb{R}^n$, and any discrete-time input $u_d:\llbracket 0,\ell \rrbracket \rightarrow \mathbb{R}^m$, there exists a continuous-time input $u_c \in \mathbb{I}_N^\tau (u_d)$ such that 
	\begin{equation}\label{InterpolatingCondition}
		\restr{x_c (\cdot; x_0,u_c)}{[0,\ell\tau]} \in \mathbb{I}_N^\tau \Big(\restr{x_d(\cdot;x_0,u_d)}{\llbracket 0,\ell\rrbracket}\Big).
	\end{equation}
\end{definition}
\vspace*{2mm}

It follows from Definition~\ref{DefInterpolatingSystem} that for any input-state trajectory in $\bm{\Sigma}_d$, there exists one in $\bm{\Sigma}_c$ that is its $N$-th order interpolation with respect to $\tau$. In other words, when $\bm{\Sigma}_c$ is an $N$-th order interpolator of $\bm{\Sigma}_d$, it is able to generate an $N$-th order interpolation of any input-state trajectory of $\bm{\Sigma}_d$. 

We now employ the shifted Legendre polynomials to obtain an \emph{algebraic} characterization of system interpolation. Specifically, by making use of the polynomial functions \eqref{PhiFunction}, we construct the vectors
\begin{subequations}\label{Necessarymatrices}
	\begin{equation}
		\begin{aligned}
			\Phi \!= \!\begin{bmatrix}
				\phi_1(\tau)\\\phi_2(\tau)\\\vdots\\\phi_N(\tau)
			\end{bmatrix}, \ \psi \!=\! \begin{bmatrix}
				\dot{\phi}_1(t_{N,0})\\\dot{\phi}_2(t_{N,0})\\\vdots\\\dot{\phi}_{N}(t_{N,0})\\
			\end{bmatrix}, \ \sigma \!=\! \begin{bmatrix}
				\dot{\phi}_0(t_{N,1})\\\dot{\phi}_0(t_{N,2})\\\vdots\\\dot{\phi}_0(t_{N,N})
			\end{bmatrix}, 
		\end{aligned}
	\end{equation}
	\text{and the matrix}
	\begin{equation}
		\Psi \!= \!\begin{bmatrix}
			\dot{\phi}_1(t_{N,1})&\dot{\phi}_1(t_{N,2})&\cdots&\dot{\phi}_1(t_{N,N})\\
			\dot{\phi}_2(t_{N,1})&\dot{\phi}_2(t_{N,2})&\cdots&\dot{\phi}_2(t_{N,N})\\
			\vdots&\vdots&\ddots&\vdots\\
			\dot{\phi}_1(t_{N,1})&\dot{\phi}_1(t_{N,2})&\cdots&\dot{\phi}_N(t_{N,N})
		\end{bmatrix},
	\end{equation}
\end{subequations}
which facilitate the characterization of system interpolation as a subspace inclusion (see \cite[Theorem 1]{Pirastzehzad2025} for further details). 
\begin{proposition}\label{PropCharacterization}
   $\bm{\Sigma}_c$ is an $N$-th order interpolator of $\bm{\Sigma}_d$ with respect to $\tau$ if and only if 
	\begin{equation}\label{Finalcharacterization}
		\begin{aligned}
			\operatorname{Im}\!\!\begingroup 
			\setlength\arraycolsep{1.25pt}\begin{bmatrix}
				A_c-\dot{\phi}_0(0)I&B_c&0\\
				-\sigma\otimes I&0&0\\
				A_d - \phi_0(\tau)I&B_d&0\\
				0&\!\!\!-\phi_0(\tau)I&I
			\end{bmatrix}\endgroup \!\subset
			\operatorname{Im}\!\!\begingroup 
			\setlength\arraycolsep{2pt}\begin{bmatrix}
				\psi^\top \otimes I & 0\\
				\Psi^\top\oplus(-A_c)&-I \otimes B_c\\
				\Phi^\top \otimes I&0\\
				0 &\Phi^\top \otimes I
			\end{bmatrix} 
			\endgroup\!,
		\end{aligned}
	\end{equation}
	where $\phi_0$ is given by \eqref{PhiFunction} and $\Phi$, $\psi$, $\sigma$, and $\Psi$ are defined in  \eqref{Necessarymatrices}.  
\end{proposition}
\vspace*{2mm}

Proposition~\ref{PropCharacterization} also enables us to characterize the set of \emph{all} continuous-time input-state trajectories that establish \eqref{InterpolatingCondition}. We thus suppose that $\bm{\Sigma}_c$ is an $N$-th order interpolator of $\bm{\Sigma}_d$ with respect to $\tau$. Given any integer $\ell\geq 1$, any initial condition $x_0\in\mathbb{R}^n$, and any control sequence $u_d:\llbracket 0,\ell\rrbracket\rightarrow\mathbb{R}^m$, for any integer $0\leq i\leq \ell-1$, we let $X^i\in\mathbb{R}^{n\times N}$ and $U^i\in\mathbb{R}^{m\times N}$ be such that 
\begin{equation}\label{AlgorithmTechnical1}
	\begin{aligned}
		&\begingroup 
		\setlength\arraycolsep{2pt}\begin{bmatrix}
			\psi^\top \otimes I & 0\\
			\Psi^\top\oplus (-A_c)&-I \otimes B_c\\
			\Phi^\top \otimes I&0\\
			0 &\Phi^\top \otimes I
		\end{bmatrix} 
		\endgroup \begin{bmatrix}
			\operatorname{vec}(X^i)\\\operatorname{vec}(U^i)
		\end{bmatrix}\\
		& = \begingroup 
		\setlength\arraycolsep{2pt}\begin{bmatrix}
			A_c-\dot{\phi}_0(0)I&B_c&0\\
			-\sigma\otimes I&0&0\\
			A_d - \phi_0(\tau)I&B_d&0\\
			0&-\phi_0(\tau)I&I
		\end{bmatrix}\endgroup \begin{bmatrix}
			x_d(i;x_0,u_d)\\u_d(i)\\u_d(i+1)
		\end{bmatrix}.
	\end{aligned}
\end{equation}
We note from Proposition~\ref{PropCharacterization} that such $X^i$ and $U^i$ always exist, while they are not necessarily unique. After partitioning $X^i$ and $U^i$ as $X^i = [X_1^i\,X_2^i\,\cdots\,X_N^i]$ and $U^i = [U_1^i\,U_2^i\,\cdots\,U_N^i]$, respectively, we define signals $\bar{x}_c^i\in\restr{\mathbb{P}_N^n}{[0,\tau]}$ and $\bar{u}_c^i\in\restr{\mathbb{P}_N^m}{[0,\tau]}$ as
\begin{equation}\label{AlgorithmTechnical2}
 \begin{aligned}
 	\bar{u}_c^i(t) &= \phi_0(t) u_d(i) + \sum_{j=1}^{N}\phi_j(t)U_j^i,\\
 	\bar{x}_c^i(t) &= \phi_0(t) x_d(i;x_0,u_d) + \sum_{j=1}^{N}\phi_j(t)X_j^i, 
 \end{aligned}
\end{equation}
for all $t\in[0,\tau]$. We then recall \eqref{SpectralRepresentation} and observe from \eqref{AlgorithmTechnical2} that  
\begin{equation*}
	\begin{aligned}
		\bar{u}_c^i(0) &= u_d(i), & \bar{u}_c^i(\tau) &= u_d(i+1),\\
		\bar{x}_c^i(0) &= x_d(i;x_0,u_d), & \bar{x}_c^i(\tau) &= x_d(i+1;x_0,u_d).
	\end{aligned}
\end{equation*}
It follows from Proposition~\ref{PropCharacterization} that the pair $(\bar{u}_c^i,\bar{x}_c^i)$ is a polynomial solution to $\bm{\Sigma}_c$ over the subinterval $[i\tau,(i+1)\tau]$. We collect all such solutions in the set 
	\begin{align}
		\mathcal{T}_c^i(x_0,u_d) = \Big\{(\bar{u}_c^i,\bar{x}_c^i) &\left\vert \bar{u}_c^i\in\restr{\mathbb{P}_N^m}{[0,\tau]},  \bar{x}_c^i\in \restr{\mathbb{P}_N^n}{[0,\tau]},  \right. \nonumber\\
		& \;\;\text{$\exists X^i\!,U^i$ s.t. \eqref{AlgorithmTechnical1}--\eqref{AlgorithmTechnical2} hold} \Big\}.  \label{SetTci}
	\end{align}
We now concatenate polynomial solutions from successive subintervals to construct piecewise polynomial solutions to $\bm{\Sigma}_c$ over the interval $[0,\ell\tau]$. We therefore let $(\bar{u}_c^i,\bar{x}_c^i)\in\mathcal{T}_c^i(x_0,u_d)$ and define the continuous-time signals $\bar{u}_c:[0,\ell \tau]\rightarrow\mathbb{R}^m$ and $\bar{x}_c:[0,\ell \tau]\rightarrow\mathbb{R}^n$ such that 
\begin{equation}\label{Stepwise}
	\begin{aligned}
		\forall t \in [0,\tau]: \;  \bar{u}_c(i\tau + t) = \bar{u}_c^i(t), \; \bar{x}_c(i\tau + t) = \bar{x}_c^i(t),
	\end{aligned}
\end{equation}
for all $i=0,1,\ldots,\ell-1$. We observe that the pair $(\bar{u}_c,\bar{x}_c)$ is a piecewise polynomial solution to $\bm{\Sigma}_c$ over the interval $[0,\ell\tau]$. After collecting these solutions into the set
\begin{align}
	\mathcal{T}_c(x_0,u_d) = \Big\{(\bar{u}_c,\bar{x}_c) \big\vert & \bar{u}_c:[0,\ell \tau]\rightarrow\mathbb{R}^m,\bar{x}_c:[0,\ell \tau]\rightarrow\mathbb{R}^n,\nonumber\\
	&\forall i\in\llbracket0,\ell-1\rrbracket, \text{$\exists (\bar{x}_c^i,\bar{u}_c^i)\in\mathcal{T}_c^i(x_0,u_d)$}\nonumber\\
	&\text{s.t. \eqref{Stepwise} holds}\Big\},\label{SetTc}
\end{align}
we obtain the following result, which is a restatement of Theorem 2 in \cite{Pirastzehzad2025}. 
\begin{theorem}\label{ThTrajectorySetCharacterization}
	Let $\bm{\Sigma}_c$ be an $N$-th order interpolator of $\bm{\Sigma}_d$ with respect to $\tau$. Then, for any integer $\ell \geq 1$, any initial condition $x_0\in\mathbb{R}^n$, and any discrete-time input $u_d:\llbracket 0,\ell \rrbracket \rightarrow \mathbb{R}^m$, we have that
	\begin{equation*}
		\begin{aligned}
			\bigg(u_c, \restr{x_c(\cdot;x_0,u_c)}{[0,\ell \tau]}\bigg)& \in \mathcal{T}_c(x_0,u_d) \\
			&\Longleftrightarrow  \text{$u_c\in\mathbb{I}_N^\tau (u_d)$ s.t. \eqref{InterpolatingCondition} holds.}
		\end{aligned}
	\end{equation*}
\end{theorem}
\vspace*{2mm}

Theorem~\ref{ThTrajectorySetCharacterization} allows for the construction of all continuous-time input-state trajectories that establish \eqref{InterpolatingCondition}. Specifically, for any given $u_d:\llbracket 0,\ell \rrbracket \rightarrow \mathbb{R}^m$, one may apply the same procedure used in the derivation of $\mathcal{T}_c(x_0,u_d)$ to generate any such trajectory. It is in fact this procedure that we will later use to conduct continuous-time execution. 

\subsubsection{Discretization by System Interpolation}
We now discretize the continuous-time dynamics $\bm{\Sigma}_c$ according to the notion of system interpolation to obtain its discrete-time model $\bm{\Sigma}_d$. Specifically, given the continuous-time system \eqref{ContinuousSystem}, for a sampling time $\tau>0$ and an integer $N\geq 1$, we construct a discrete-time system \eqref{DiscreteSystem} in such a way that \eqref{ContinuousSystem} is an $N$-th order interpolator of \eqref{DiscreteSystem} with respect to $\tau$. 

We therefore treat $A_d$ and $B_d$ as \emph{unknown} parameters. In light of Proposition~\ref{PropCharacterization}, we seek matrices $A_d$ and $B_d$ such that \eqref{Finalcharacterization} holds. Following the procedure described in \cite[Section V]{Pirastzehzad2025}, we construct the matrices
\begin{subequations}\label{QRSV}
	\begin{equation}
		\begin{aligned}
			Q_N^\tau (A_c,B_c) &= \begin{bmatrix}
				\psi^\top \otimes I & 0\\
				\Psi^\top\oplus(-A_c)&-I \otimes B_c\\
				\Phi^\top \otimes I&0\\
				0 &\Phi^\top \otimes I
			\end{bmatrix},\\
			R_N^\tau(A_c,B_c) &= \begingroup\setlength\arraycolsep{1.25pt}\begin{bmatrix}
				A_c-\dot{\phi}_0(0)I&B_c&0\\
				-\sigma\otimes I&0&0\\
				- \phi_0(\tau)I&0&0\\
				0&\!\!\!-\phi_0(\tau)I&I
			\end{bmatrix}\endgroup,
		\end{aligned}
	\end{equation}
	\text{and}
	\begin{equation}
		\begin{aligned}
			T_1 &= \begin{bmatrix}
				0\\0\\I\\0
			\end{bmatrix}, &T_2&= \begin{bmatrix}
				I&0&0\\
				0&I&0
			\end{bmatrix},
		\end{aligned}
	\end{equation}
\end{subequations}
which enable us to equivalently express \eqref{Finalcharacterization} as the existence of a matrix $M$ such that 
\begin{equation}\label{CompactRep}
	Q_N^\tau (A_c,B_c) M = R_N^\tau(A_c,B_c) + T_1\begin{bmatrix}
		A_d&B_d
	\end{bmatrix}T_2.
\end{equation}
We emphasize that the matrices $Q_N^\tau (A_c,B_c)$ and $R_N^\tau(A_c,B_c)$ are completely in terms of known parameters (\textit{i.e.,} the matrices $A_c,B_c$, the polynomial degree $N$, and the sampling time $\tau$). After  defining the vectors
\begin{equation}\label{IntermediateResult}
	\begin{aligned}
		v_1 &= \operatorname{vec}(M), &v_2&=\operatorname{vec}\Big(\begin{bmatrix}
			A_d&B_d
		\end{bmatrix}\Big),
	\end{aligned}
\end{equation} 
we vectorize \eqref{CompactRep} to obtain the following result, which is taken from \cite[Theorem 3]{Pirastzehzad2025}. 
\begin{theorem}\label{Th_Design}
	Given a continuous-time system \eqref{ContinuousSystem}, there exist matrices $A_d$ and $B_d$ such that \eqref{ContinuousSystem} is an $N$-th order interpolator of \eqref{DiscreteSystem} with respect to $\tau$ if and only if there exist vectors $v_1$ and $v_2$ such that 
	\begin{equation}\label{DesignResult}
		\begin{bmatrix}
			I\otimes Q_N^\tau (A_c,B_c) &-T_2^\top \otimes T_1
		\end{bmatrix}\begin{bmatrix}
			v_1\\
			v_2\\
		\end{bmatrix} = \operatorname{vec}\big( R_N^\tau(A_c,B_c) \big).
	\end{equation}
	Moreover, for any vectors $v_1$ and $v_2$ that satisfy \eqref{DesignResult}, the matrices $A_d$ and $B_d$ recovered from \eqref{IntermediateResult} are such that \eqref{ContinuousSystem} is an $N$-th order interpolator of \eqref{DiscreteSystem} with respect to $\tau$.
\end{theorem}
\vspace*{2mm}

Theorem~\ref{Th_Design} characterizes the existence of $\bm{\Sigma}_d$ as solvability of the linear equation \eqref{DesignResult}, whose solution gives the parameters of $\bm{\Sigma}_d$, \textit{i.e.,} matrices $A_d$ and $B_d$. 

We may now perform the \emph{discretization} step of the proposed hierarchical control scheme in accordance with the methodology discussed above. Given the continuous-time system $\bm{\Sigma}_c$ and its corresponding STL formula $\mathfrak{S}_c$, we first select the sampling time $\tau>0$ such that every time bound specified in $\mathfrak{S}_c$ appears as an integer multiple of $\tau$. We then apply the recursive procedure \eqref{SpecificationDiscretizationProcedure} to construct the discrete-time STL formula $\mathfrak{S}_d$.  Lastly, for a polynomial degree $N\geq 1$, we utilize Theorem~\ref{Th_Design} to construct a discrete-time model $\bm{\Sigma}_d$ such that $\bm{\Sigma}_c$ is an $N$-th order interpolator of $\bm{\Sigma}_d$ with respect to $\tau$. 
\section{Discretization-aware Planning}\label{Sec_Planning}
Having conducted discretization, we now proceed with the \emph{planning} step of the hierarchical control scheme. We thus suppose that the discrete-time STL formula $\mathfrak{S}_d$ is obtained according to the recursive procedure \eqref{SpecificationDiscretizationProcedure} and that the discrete-time system $\bm{\Sigma}_d$ is constructed such that $\bm{\Sigma}_c$ is an $N$-th order interpolator of $\bm{\Sigma}_d$ with respect to $\tau$. 

For an initial condition $x_0\in\mathbb{R}^n$, our objective is therefore to compute a discrete-time control sequence $u_d:\llbracket0,\mathfrak{h}_{\mathfrak{S}_d}\rrbracket\rightarrow\mathcal{U}$ that establishes \eqref{DiscreteControlProblem}. Importantly, this control sequence must be computed in such a way that the continuous-time control input $u_c$, which is obtained via continuous-time execution, establishes \eqref{ControlProblem} while satisfying the input constraints (\textit{i.e.,} $u_c(t) \in \mathcal{U}$ for all $t$). We guarantee this by imposing additional constraints on the control sequence $u_d$. 

We split the derivation of these constraints into two parts. In the first part, by exploiting the fact that input-state trajectories in $\bm{\Sigma}_c$ can be realized as piecewise polynomial interpolations of input-state trajectories in $\bm{\Sigma}_d$, we formulate constraints that enable us to transform a control sequence $u_d:\llbracket0,\mathfrak{h}_{\mathfrak{S}_d}\rrbracket\rightarrow\mathcal{U}$ that yields \eqref{DiscreteControlProblem} into a continuous-time control input $u_c:[0,\mathfrak{h}_{\mathfrak{S}_c}]\rightarrow\mathcal{U}$ that establishes \eqref{ControlProblem}. In the second part, on the other hand, we adopt the MILP encoding proposed in \cite{7039363} to formulate \eqref{DiscreteControlProblem} as mixed-integer linear constraints in terms of the initial condition $x_0$ and the control sequence $u_d$. 

\subsection{Control Constraints}\label{Sec_Constraints}
We derive constraints on a control sequence $u_d$ such that 1) its piecewise polynomial interpolation yields a continuous-time control input $u_c$ that takes values in the set $\mathcal{U}$, and 2) whenever this control sequence $u_d$ satisfies \eqref{DiscreteControlProblem}, then the corresponding continuous-time control input $u_c$ establishes \eqref{ControlProblem}. 

We begin by deriving a condition on $u_d$ that guarantees $u_c$ remains in the set $\mathcal{U}$. Given an initial condition $x_0\in\mathbb{R}^n$, for an integer $\ell\geq 1$, we consider a control sequence $u_d:\llbracket0,\ell\rrbracket\rightarrow\mathcal{U}$ and let the continuous-time input $u_c\in\mathbb{I}_N^\tau(u_d)$ be such that \eqref{InterpolatingCondition}. We thus conclude from Theorem~\ref{ThTrajectorySetCharacterization} that
\begin{equation}\label{KeyInclusion}
	\bigg(u_c,\restr{x_c(\cdot;x_0,u_c)}{[0,\ell \tau]}\bigg) \in \mathcal{T}_c(x_0,u_d),
\end{equation}
which, as a result of \eqref{SetTc}, implies that for all $i\in\llbracket0,\ell-1\rrbracket$, there exist matrices  $X^i\!\in\mathbb{R}^{n\times N}$ and $U^i\!\in\mathbb{R}^{m\times N}$ such that \eqref{AlgorithmTechnical1} holds, while for all $t\in[0,\tau]$, we have that
\begin{subequations}\label{SpectralRep}
	\begin{align}
		u_c(i\tau+t) \!&=\! \phi_0(t) u_d(i) + \sum_{j=1}^{N}\phi_j(t)U_j^i,\label{SpectralRepU}\\
		x_c(i\tau+t; x_0,u_c) \!&=\! \phi_0(t) x_d(i;x_0,u_d) + \sum_{j=1}^{N}\phi_j(t)X_j^i, \label{SpectralRepX}
	\end{align}
\end{subequations}
where $X^i = [X_1^i\,X_2^i\,\cdots\,X_N^i]$ and $U^i = [U_1^i\,U_2^i\,\cdots\,U_N^i]$. We will therefore derive a condition on $u_d$ guaranteeing that
\begin{equation}\label{InputConstraint}
	\forall t\in[0,\tau]: u_c(i\tau+t) \in \mathcal{U},
\end{equation}
for all $i\in\llbracket0,\ell-1\rrbracket$. We accomplish this by formulating \eqref{InputConstraint} in terms of the so-called \emph{Bernstein basis polynomials} \cite{lorentz2012bernstein}, which are defined as
\begin{equation*}
	\mathfrak{b}_{k,N}(t) = \frac{N!}{k!(N-k)!}t^k(1-t)^{N-k}, \quad k=0,1,\ldots,N.
\end{equation*}
These polynomials are non-negative over the interval $[0,1]$ (\textit{i.e.,} $\mathfrak{b}_{k,N}(t) \geq 0$ for all $t\in[0,1]$) and form a partition of unity (\textit{i.e.,} $\sum_{k=0}^{N}\mathfrak{b}_{k,N}(t) = 1$). Moreover, they are linearly independent (over the interval $[0,1]$) and, therefore, form a basis for the space of polynomials up to degree $N$ (over the interval $[0,1]$) \cite{cargo1966bernstein}. This, due to the fact that $\restr{u_c}{[i\tau,(i+1)\tau]}\in\restr{\mathbb{P}_N^m}{[0,\tau]}$, implies the existence of vectors $v_0^i,v_1^i,\ldots,v_N^i\in\mathbb{R}^m$ such that 
\begin{equation}\label{InputBernsteinRep}
	\forall t\in[0,\tau]: \quad u_c(i\tau+t) = \sum_{k=0}^{N} v_k^i \mathfrak{b}_{k,N}\left(\frac{t}{\tau}\right).
\end{equation}
Recalling that the polynomials $\mathfrak{b}_{k,N}$ are non-negative and form a partition of unity, we observe from \eqref{InputBernsteinRep} that $u_c(i\tau+t)$ can be written as a \emph{convex} combination of the vectors $v_0^i,v_1^i,\ldots,v_N^i$. It accordingly follows from the convexity of $\mathcal{U}$ that in order to establish \eqref{InputConstraint}, it suffices to ensure that $v_k^i\in\mathcal{U}$ for all $k=0,1,\ldots,N$. We write the vectors $v_k^i$ in terms of the matrix $U^i$ to obtain a condition on $u_d$ that guarantees \eqref{InputConstraint}. 

Bearing \eqref{Dirac} in mind, we observe from \eqref{SpectralRepU} that $u_c(i\tau + t_{N,0}) = u_d(i)$ and that $u_c(i\tau+t_{N,j}) = U_j^i$ for $j=1,2,\ldots,N$. We then utilize \eqref{InputBernsteinRep} to conclude that
\begin{subequations}\label{UEquations}
	\begin{equation}
		\sum_{k=0}^{N}v_k^i\mathfrak{b}_{k,N}\left(\frac{t_{N,0}}{\tau}\right) = u_d(i),
	\end{equation}
	\text{and that}
	\begin{equation}
		\sum_{k=0}^{N}v_k^i\mathfrak{b}_{k,N}\left(\frac{t_{N,j}}{\tau}\right) = U_j^i,
	\end{equation}
\end{subequations}
for all $j=1,2,\ldots,N$. We accordingly construct the matrix
\begin{equation}\label{B}
	\mathfrak{B} = \begingroup\setlength\arraycolsep{1.25pt}\begin{bmatrix}
		\mathfrak{b}_{0,N}\left(\frac{t_{N,0}}{\tau}\right)&\mathfrak{b}_{0,N}\left(\frac{t_{N,1}}{\tau}\right)&\cdots&\mathfrak{b}_{0,N}\left(\frac{t_{N,N}}{\tau}\right)
		\vspace*{0.5mm}\\
		\mathfrak{b}_{1,N}\left(\frac{t_{N,0}}{\tau}\right)&\mathfrak{b}_{1,N}\left(\frac{t_{N,1}}{\tau}\right)&\cdots&\mathfrak{b}_{1,N}\left(\frac{t_{N,N}}{\tau}\right)\\
		\vdots&\vdots&\ddots&\vdots\\
		\mathfrak{b}_{N,N}\left(\frac{t_{N,0}}{\tau}\right)&\mathfrak{b}_{N,N}\left(\frac{t_{N,1}}{\tau}\right)&\cdots&\mathfrak{b}_{N,N}\left(\frac{t_{N,N}}{\tau}\right)\\
	\end{bmatrix}\endgroup,
\end{equation}
and collect the equations \eqref{UEquations} into
\begin{equation*}
	\begin{bmatrix}
		v_0^i&v_1^i&\cdots&v_N^i
	\end{bmatrix}\mathfrak{B} = \begin{bmatrix}
		u_d(i)&U^i
	\end{bmatrix}.
\end{equation*}
Noting from Lemma~\ref{TechnicalLemma} in the appendix that $\mathfrak{B}$ is invertible, we obtain
\begin{equation}
	v_k^i = \begin{bmatrix}
		u_d(i)&U^i
	\end{bmatrix} \mathfrak{B}^{-1}e_k,
\end{equation}
where $e_k$ is the standard basis vector whose entries are all 0 except for the $k$-th one, which is 1. We therefore guarantee \eqref{InputConstraint} by establishing that 
\begin{equation}\label{InputGuarantee}
	\forall k\in\llbracket0,N\rrbracket: \quad \begin{bmatrix}
		u_d(i)&U^i
	\end{bmatrix} \mathfrak{B}^{-1}e_k \in \mathcal{U}.
\end{equation}
This leads to the following result, which gives a sufficient condition that guarantees $u_c$ remains in $\mathcal{U}$. The proof of this result directly follows from the procedure described above and is therefore omitted. 
\begin{theorem}\label{ThControlConstraints}
	Let $\bm{\Sigma}_c$ be an $N$-th order interpolator of $\bm{\Sigma}_d$ with respect to $\tau$. Given an initial condition $x_0\in\mathbb{R}^n$, let the control sequence $u_d:\llbracket0,\ell\rrbracket\rightarrow\mathcal{U}$ (with $\ell\geq 1$) be such that for every $i\in\llbracket0,\ell-1\rrbracket$, there exist matrices $X^i\in\mathbb{R}^{n\times N}$ and $U^i\in\mathbb{R}^{m\times N}$ that satisfy \eqref{AlgorithmTechnical1} and \eqref{InputGuarantee}. Then, for any such $U^0, U^1, \ldots, U^{\ell-1}$, the continuous-time control input $u_c$ defined through \eqref{SpectralRep} satisfies
	 \begin{equation}\label{DefInputConstraint}
	 	\forall t\in[0,\ell\tau]: \quad u_c(t) \in \mathcal{U}.
	 \end{equation}
\end{theorem}
\vspace*{2mm}

Theorem~\ref{ThControlConstraints} provides constraints on the control sequence $u_d$ which guarantee that its piecewise polynomial interpolation produces a continuous-time input $u_c\in\mathbb{I}_N^\tau(u_d)$ that establishes \eqref{InterpolatingCondition}, while remaining in the set $\mathcal{U}$. 

We next derive a condition which ensures that piecewise polynomial interpolation of discrete-time input-state trajectories satisfying \eqref{DiscreteControlProblem} yields continuous-time ones that establish \eqref{ControlProblem}. Since the predicate \eqref{Predicate} is a fundamental constituent of the STL formula $\mathfrak{S}_c$, we obtain such a condition by characterizing the satisfaction of \eqref{Predicate} by a continuous-time input-state trajectory in terms of its discrete-time counterpart. 

With $x_0$, $u_d$, and $u_c$ defined as above, let a predicate \eqref{Predicate} be such that 
\begin{equation}\label{MidCond}
	\Big(\restr{x_d(\cdot;x_0,u_d)}{\llbracket0,\ell\rrbracket},i\Big) \models \mu,
\end{equation}
for some $i\in\llbracket0,\ell-1\rrbracket$. We will derive a condition on $x_0$ and $u_d$ which guarantees that 
	\begin{equation}\label{StateConstraint}
		\forall t\in[0,\tau]: \quad \Big(\restr{x_c(\cdot;x_0,u_c)}{[0,\ell\tau]},i\tau+t\Big) \models \mu,
	\end{equation}
indicating that the continuous-time state trajectory satisfies \eqref{Predicate} throughout the interval $[i\tau,(i+1)\tau]$. To accomplish this, we define the scalar function $P^{i}:[0,\tau]\rightarrow\mathbb{R}$ such that 
\begin{equation}\label{P}
	\begin{aligned}
		\forall t\in[0,\tau]:\;	P^{i} (t) \!= \!\gamma^\top\!\!\left(\restr{x_c(i\tau+t;x_0,u_c)}{[0,\tau]}\right)+\beta.
	\end{aligned}
\end{equation}
We then observe from the semantics \eqref{Semantics} that \eqref{StateConstraint} is equivalent to
\begin{equation}\label{PolyCond}
	\min_{t\in[0,\tau]} P^{i}(t)>0.
\end{equation}
To establish this, we note from \eqref{SpectralRepX} that $P^{i}\in\restr{\mathbb{P}_N^1}{[0,\tau]}$. There consequently exist constants $c_0^{i},c_1^{i},\ldots,c_N^{i}\in\mathbb{R}$ such that 
\begin{equation}\label{BernsteinExpansion}
	\forall t\in[0,\tau]: \quad P^{i}(t) = \sum_{k=0}^{N} c_k^{i}\mathfrak{b}_{k,N}\left(\frac{t}{\tau}\right).
\end{equation}
Bearing in mind that the polynomials $\mathfrak{b}_{k,N}$ are non-negative and form a partition of unity, we conclude from \eqref{BernsteinExpansion} that 
\begin{equation}\label{MinCond}
	\min_{t\in[0,\tau]}P^{i}(t) \geq \min_{k\in\llbracket0,N\rrbracket} c_k^{i}.
\end{equation}
We compute the constants $c_k^{i}$ in terms of the matrix $X^{i}$. After substituting \eqref{SpectralRepX} in \eqref{P}, we utilize \eqref{Dirac} and \eqref{BernsteinExpansion} to conclude that 
	\begin{align*}
		\sum_{k=0}^{N}c_k^{i}\mathfrak{b}_{k,N}\left(\frac{t_{N,0}}{\tau}\right) &= \gamma^\top x_d(i;x_0,u_d) + \beta,
	\intertext{and}
		\sum_{k=0}^{N}c_k^{i}\mathfrak{b}_{k,N}\left(\frac{t_{N,j}}{\tau}\right) &= \gamma^\top X_j^{i} + \beta,
	\end{align*}
for all $j = 1,2,\ldots,N$. We then collect these equations into 
\begin{equation*}
	\begin{bmatrix}
		c_0^{i}&c_1^{i}&\cdots&c_N^{i}
	\end{bmatrix}\mathfrak{B} = \gamma^\top\!\!\begin{bmatrix}
	x_d(i;x_0,u_d)&
	X^{i}
	\end{bmatrix} + \beta\bm{1}^\top,
\end{equation*} 
where $\bm{1}$ is a vector of ones with appropriate dimension. This immediately gives 
\begin{equation*}
	c_k^{i} = \Big(\gamma^\top\!\!\begin{bmatrix}
		x_d(i;x_0,u_d)&
		X^{i}
	\end{bmatrix} + \beta\bm{1}\Big) \mathfrak{B}^{-1}e_k,
\end{equation*}
which, together with \eqref{PolyCond} and \eqref{MinCond}, implies that \eqref{StateConstraint} holds if
\begin{equation}\label{GuaranteeCondition}
	\forall k\in\llbracket0,N\rrbracket: \Big(\!\gamma^\top\!\!\begin{bmatrix}
		x_d(i;x_0,u_d)&\!\!\! X^{i}
	\end{bmatrix} + \beta\bm{1}\!\Big) \mathfrak{B}^{-1}e_k>0.
\end{equation}
This leads to the following result, which characterizes a condition under which satisfaction of a predicate \eqref{Predicate} by a discrete-time input-state trajectory at an instant $i$ ensures its satisfaction by the corresponding continuous-time input-state trajectory over the interval $[i\tau,(i+1)\tau]$. The proof of this result is an immediate consequence of the preceding derivation and is therefore omitted.
\begin{theorem}\label{ThGuaranteeCondition}
	Let $\bm{\Sigma}_c$ be an $N$-th order interpolator of $\bm{\Sigma}_d$ with respect to $\tau$. Given an initial condition $x_0\in\mathbb{R}^n$ and a predicate \eqref{Predicate}, let the control sequence $u_d: \llbracket0,\ell\rrbracket\rightarrow \mathcal{U}$ (with $\ell\geq1$) be such that \eqref{MidCond} holds for some $i\in\llbracket0,\ell-1\rrbracket$ and that there exist matrices $X^{i}\in\mathbb{R}^{n\times N}$ and $U^{i}\in\mathbb{R}^{m\times N}$ that satisfy \eqref{AlgorithmTechnical1} and \eqref{GuaranteeCondition}. Then, for any such $U^0, U^1,\ldots,U^{\ell-1}$, the continuous-time control input $u_c$ defined through \eqref{SpectralRep} establishes \eqref{StateConstraint}. 
\end{theorem}
\vspace*{2mm}

Theorem~\ref{ThGuaranteeCondition} provides constraints on $x_0$ and $u_d$ which guarantee that whenever $u_d$ enforces a predicate \eqref{Predicate} at an instant $i$, its piecewise polynomial interpolation produces a continuous-time input $u_c\in\mathbb{I}_N^\tau(u_d)$ that establishes \eqref{InterpolatingCondition}, while ensuring that the predicate \eqref{Predicate} is satisfied over the interval $[i\tau,(i+1)\tau]$. 

We will later use Theorems~\ref{ThControlConstraints} and \ref{ThGuaranteeCondition} to formulate an MILP whose solution gives a control sequence which, after continuous-time execution, yields a control input that establishes \eqref{ControlProblem} while respecting the input constraints (\textit{i.e.,} taking values in $\mathcal{U}$). 
\subsection{STL Constraints}\label{Sec_Encoding}
We adopt the methodology proposed in \cite[Section IV]{7039363} to encode the discrete-time STL formula $\mathfrak{S}_d$ in terms of mixed-integer linear constraints. Specifically, for any $i\in\llbracket0,\mathfrak{h}_{\mathfrak{S}_d}\rrbracket$, we introduce a \emph{binary} variable $z_{\mathfrak{S}_d}^i \in \{0,1\}$ such that 
\begin{equation}\label{FundamentalSTLConstraint}
	 z_{\mathfrak{S}_d}^i = 1 \Longleftrightarrow \big(x_d(\cdot;x_0,u_d),i\big) \models \mathfrak{S}_d .
\end{equation}

We let $\mathfrak{S}$ be a subformula in $\mathfrak{S}_d$ and define the binary variable $z_{\mathfrak{S}}^i\in\{0,1\}$ that captures its validity at an instant $i$. As for a predicate \eqref{Predicate}, we take $\mathfrak{S} = \mu$ and define $z_{\mathfrak{S}}^i$ such that 
\begin{equation}\label{PredicateConstraints}
	\begin{aligned}
		\gamma^\top x_d(i;x_0,u_d) + \beta &\leq M^i z_{\mathfrak{S}}^i - \varepsilon^i,\\
		\gamma^\top x_d(i;x_0,u_d) + \beta &\geq M^i(z_{\mathfrak{S}}^i -1) + \varepsilon^i, 
	\end{aligned}
\end{equation}
where $M^i >0$ and $\varepsilon^i>0$ are respectively sufficiently large and small numbers that are a priori chosen by the designer (for further details, see \cite[Section IV-B.1]{7039363}). We thus observe that $\big(x_d(i;x_0,u_d),i\big)\models \mu$ if and only if $z_\mu^i = 1$.

We next address the Boolean operators. Considering the negation operator, we take $\mathfrak{S}  = \neg \mathscr{S}$ (for a given subformula $\mathscr{S}$ with binary variable $z_{\mathscr{S}}^i\in\{0,1\}$) and define 
\begin{equation}\label{NegationConstraint}
	z_{\mathfrak{S}}^i = 1- z_{\mathscr{S}}^i.
\end{equation}
Regarding conjunction, on the other hand, we let $\mathfrak{S} =\bigwedge_{l=1}^{k}\mathscr{S}^l$ (for given subformulas $\mathscr{S}^1,\mathscr{S}^2,\ldots,\mathscr{S}^k$) and define $z_{\mathfrak{S}}^i$ such that 
\begin{equation}\label{ConjunctionConstraint}
	\begin{aligned}
		\forall l\in\llbracket1,k\rrbracket: \quad z_{\mathfrak{S}}^i &\leq z_{\mathscr{S}^l}^i,\\
		z_{\mathfrak{S}}^i &\geq 1 - k + \sum_{l=1}^{k}z_{\mathscr{S}^l}^i,
	\end{aligned}
\end{equation}
where $z_{\mathscr{S}^l}^i\in\{0,1\}$ is the binary variable associated with the subformula $\mathscr{S}^l$ for $l=1,2,\ldots,k$. We further use the fact that $\bigvee_{l=1}^{k}\mathscr{S}^l = \neg(\bigwedge_{l=1}^{k}\neg\mathscr{S}^l)$ to define the binary variable for the disjunction operator in terms of \eqref{NegationConstraint} and \eqref{ConjunctionConstraint}. 

We lastly consider temporal operators. As for the always operator, we let $\mathfrak{S} = \operatorname{G}_{\llbracket k_1,k_2\rrbracket}\mathscr{S}$ (for a given subformula $\mathscr{S}$ and integers $0\leq k_1\leq k_2$) and define 
\begin{equation}\label{AlwaysConstraint}
	z_{\mathfrak{S}}^i = \bigwedge_{l=\min\{i+k_1,\mathfrak{h}_{\mathfrak{S}_d}\}}^{\min\{i+k_2,\mathfrak{h}_{\mathfrak{S}_d}\}} z_{\mathscr{S}}^l,
\end{equation}
where $z_{\mathscr{S}}^l\in\{0,1\}$ is the binary variable associated with $\mathscr{S}$ at instant $l$. Regarding the until operator, we let $\mathfrak{S} = \mathscr{S}^1\operatorname{U}_{\llbracket k_1,k_2\rrbracket}\mathscr{S}^2$ (for given subformulas $\mathscr{S}^1$ and $\mathscr{S}^2$) and define 
\begin{equation}\label{UntilConstraint}
	z_{\mathfrak{S}}^i = \bigvee_{l=\min\{i+k_1,\mathfrak{h}_{\mathfrak{S}_d}\}}^{\min\{i+k_2,\mathfrak{h}_{\mathfrak{S}_d}\}} \Big(z_{\mathscr{S}^2}^l \wedge \bigwedge_{l'=i}^{l'=l}z_{\mathscr{S}^1}^{l'}\Big),
\end{equation}
where $z_{\mathscr{S}^1}^{l'},z_{\mathscr{S}^2}^l\in\{0,1\}$ are binary variables that are respectively associated with the subformulas $\mathscr{S}^1$ (at instant $l'$) and $\mathscr{S}^2$ (at instant $l$). Considering that $\operatorname{F}_{\llbracket k_1,k_2\rrbracket}\mathscr{S} = \mathrm{True} \operatorname{U}_{\llbracket k_1,k_2\rrbracket} \mathscr{S}$, we may easily define the binary variable for the eventually operator on the basis of \eqref{UntilConstraint}. 

Based on these variables, we utilize the syntax \eqref{Syntax} to recursively obtain the binary variable $z_{\mathfrak{S}_d}^i$ according to \eqref{PredicateConstraints}--\eqref{UntilConstraint}. We note from induction on the structure of $\mathfrak{S}_d$ that \eqref{FundamentalSTLConstraint} holds (see, \textit{e.g.,} \cite{7039363,annurev:/content/journals/10.1146/annurev-control-053018-023717}), which, by definition, implies that \eqref{DiscreteControlProblem} holds if and only if 
\begin{equation}\label{SatisfactionConstraint}
	z_{\mathfrak{S}_d}^0 = 1.
\end{equation}
We collect all binary variables introduced by this recursive construction into a vector $z$. 

Furthermore, since the control sequence $u_d$ is restricted to take values in the set $\mathcal{U}$, it must satisfy the constraints
\begin{equation}\label{DiscreteInputConstraint}
	\forall i\in\llbracket0,\mathfrak{h}_{\mathfrak{S}_d}\rrbracket: \quad u_d(i) \in \mathcal{U},
\end{equation}
which can be represented by linear inequalities, as $\mathcal{U}$ is a convex polyhedron.

We combine \eqref{DiscreteInputConstraint} with the constraints defined by the binary vector $z$ and denote the resulting set of constraints by $\mathcal{Z}_{\mathfrak{S}_d}(z,u_d)$. A feasible solution to these constraints gives the control sequence $u_d$ that establishes \eqref{DiscreteControlProblem}. This observation is formalized in the following proposition, which is taken from \cite[Theorem 1]{annurev:/content/journals/10.1146/annurev-control-053018-023717}. 
\begin{proposition}\label{PropMILPEncoding}
	Given a discrete-time STL formula $\mathfrak{S}_d$, for an initial condition $x_0\in\mathbb{R}^n$, there exists a control sequence $u_d:\llbracket0,\mathfrak{h}_{\mathfrak{S}_d}\rrbracket:\rightarrow\mathcal{U}$ that establishes \eqref{DiscreteControlProblem} if and only if the constraints $\mathcal{Z}_{\mathfrak{S}_d}(z,u_d)$ are feasible.
\end{proposition}

In the next section, we combine the constraints $\mathcal{Z}_{\mathfrak{S}_d}(z,u_d)$ with those formulated in Theorem~\ref{ThControlConstraints} and \ref{ThGuaranteeCondition} to derive the MILP whose solution enables construction of the continuous-time control input that establishes \eqref{ControlProblem}.
\subsection{MILP}
The MILP is composed of constraints associated with the matrices $X^i\in\mathbb{R}^{n\times N}$ and $U^i\in\mathbb{R}^{m\times N}$, which are defined through the set $\mathcal{T}(x_0,u_d)$, and the binary variables introduced in $\mathcal{Z}_{\mathfrak{S}_d}(z,u_d)$.

We begin by formulating the constraints involving $X^i$ and $U^i$. For each $i\in\llbracket0,\mathfrak{h}_{\mathfrak{S}_d}-1\rrbracket$, we define the equation \eqref{AlgorithmTechnical1}. These equations enable the construction of continuous-time input-state trajectories as piecewise polynomial interpolations of the discrete-time ones. Additionally, we introduce the inclusion \eqref{InputGuarantee}, which ensures that piecewise polynomial interpolation of the control sequence $u_d$ produces a continuous-time input $u_c$ that remains in the  set $\mathcal{U}$. We emphasize that since $\mathcal{U}$ is a convex polyhedron, \eqref{InputGuarantee} can be expressed as linear constraints. 

For each predicate $\mu$ in the formula $\mathfrak{S}_d$, on the other hand, we introduce the constraint
	\begin{align}\label{IntegratedConstraint}
		\forall k\in\llbracket0,N\rrbracket: \Big(\!\gamma^\top \!\!\begin{bmatrix}
			x_d(i;x_0,u_d)& \!\!\! X^{i}
		\end{bmatrix} &+ \beta\bm{1}\!\Big) \mathfrak{B}^{-1}e_k \\  
		&\geq M^i(z_{\mu}^i - 1)  +  \varepsilon^i,\nonumber
	\end{align}
where $\gamma\in\mathbb{R}^n$ and $\beta\in\mathbb{R}$ are given by \eqref{Predicate}, the matrix $X^i$ satisfies \eqref{AlgorithmTechnical1}, and the binary variable $z_{\mu}^i\in\{0,1\}$ and constants $M^i,\varepsilon^i>0$ are defined as in \eqref{PredicateConstraints}. We observe from \eqref{IntegratedConstraint} that if the predicate $\mu$ is valid at instant $i$ (\textit{i.e.,} $z_\mu^i = 1$), then \eqref{GuaranteeCondition} holds. This, as a consequence of Theorem~\ref{ThGuaranteeCondition}, ensures that whenever the predicate $\mu$ is valid at an instant $i$, it remains valid throughout the interval $[i\tau,(i+1)\tau]$. We emphasize that in contrast to \eqref{PredicateConstraints}, which is necessary and sufficient for the satisfaction of $\mu$ in discrete time, \eqref{IntegratedConstraint} provides only a sufficient condition for its satisfaction in continuous time.

 We collect the matrices $X^i$ and $U^i$, introduced by \eqref{AlgorithmTechnical1}, \eqref{InputGuarantee}, \eqref{IntegratedConstraint}, respectively, into matrices $X$ and $U$ such that 
 \begin{equation}\label{Partition}
 	\begin{aligned}
 		X = \begin{bmatrix}
 			X^0&X^1&\cdots&X^{\mathfrak{h}_{\mathfrak{S}_d}-1}
 		\end{bmatrix},\\
 		U = \begin{bmatrix}
 			U^0&U^1&\cdots&U^{\mathfrak{h}_{\mathfrak{S}_d}-1}
 		\end{bmatrix}.
 	\end{aligned}
 \end{equation}

Next, we formulate the constraints involving the binary variables. Having associated each predicate $\mu$ in the formula $\mathfrak{S}_d$ with the mixed-integer linear constraints \eqref{PredicateConstraints} and \eqref{IntegratedConstraint}, we apply the recursive procedure described in the previous section to obtain the set of constraints $\mathcal{Z}_{\mathfrak{S}_d}(z,u_d)$, which capture the validity of $\mathfrak{S}_d$ at each instant $i\in\llbracket0,\mathfrak{h}_{\mathfrak{S}_d}\rrbracket$. It is worth mentioning that the constraints \eqref{DiscreteInputConstraint} need not be introduced separately as they are already enforced by \eqref{InputGuarantee}.

We therefore formulate an MILP whose constraints are 1) the equations \eqref{AlgorithmTechnical1} defined for all $i\in\llbracket0,\mathfrak{h}_{\mathfrak{S}_d}\rrbracket$, 2) the inclusion \eqref{InputGuarantee} defined for all $i\in\llbracket0,\mathfrak{h}_{\mathfrak{S}_d}\rrbracket$, and 3) the mixed-integer linear constraints that are recursively defined according to \eqref{PredicateConstraints}--\eqref{IntegratedConstraint}. We denote this \emph{augmented} MILP by $\mathcal{Z}_{\mathfrak{S}_d}^{\mathrm{aug}}(u_d,X,U,z)$, where the superscript $\mathrm{aug}$ highlights that the constraints \eqref{PredicateConstraints}--\eqref{SatisfactionConstraint} have been extended with \eqref{AlgorithmTechnical1}, \eqref{InputGuarantee}, and \eqref{IntegratedConstraint}. We further emphasize that the \emph{unknown} variables in this MILP are the control sequence $u_d$, the matrices $X^i,U^i$ (which are respectively collected in $X$ and $U$), and the binary variables recursively defined according to \eqref{PredicateConstraints}--\eqref{IntegratedConstraint}, which are stacked into the binary vector $z$. 

In the next section, we show that a feasible solution to the MILP $\mathcal{Z}_{\mathfrak{S}_d}^{\mathrm{aug}}(u_d,X,U,z)$ enables the construction of a continuous-time control input that establishes \eqref{ControlProblem}. 
\section{Continuous-time Execution}\label{Sec_CTExecution}
Having conducted planning, we now perform the last step of the hierarchical scheme, which is to construct the continuous-time control input $u_c$ that establishes \eqref{ControlProblem}. We therefore suppose that the MILP $\mathcal{Z}_{\mathfrak{S}_d}^{\mathrm{aug}}(u_d,X,U,z)$ is feasible. We partition matrices $X$ and $U$ according to \eqref{Partition}. Considering that $\bm{\Sigma}_c$ is an $N$-th order interpolator of $\bm{\Sigma}_d$ with respect to $\tau$, we construct the continuous-time control input $u_c:[0,\mathfrak{h}_{\mathfrak{S}_c}]\rightarrow\mathbb{R}^m$ according to \eqref{SpectralRepU}. We note from \eqref{Horizons} that the construction of such control input is always possible. We then immediately conclude from Theorem~\ref{ThTrajectorySetCharacterization} that \eqref{SpectralRepX} holds. We now show that the input \eqref{SpectralRepU} remains in the set $\mathcal{U}$ and that the state trajectory \eqref{SpectralRepX} satisfies the continuous-time STL formula $\mathfrak{S}_c$.

\begin{theorem}\label{ThMain}
	Consider the continuous-time system $\bm{\Sigma}_c$ and the STL formula $\mathfrak{S}_c$ constructed according to the syntax \eqref{Syntax}. Let $\tau>0$ be a sampling time such that every time bound specified in $\mathfrak{S}_c$ appears as an integer multiple of $\tau$. Let $\mathfrak{S}_d$ denote the discrete-time STL formula obtained according to the recursive procedure \eqref{SpecificationDiscretizationProcedure}. Additionally, let $\bm{\Sigma}_d$ be a discrete-time system such that $\bm{\Sigma}_c$ is an $N$-th order interpolator of $\bm{\Sigma}_d$ with respect to $\tau$. 
	
	For an initial condition $x_0\in\mathbb{R}^n$, there exists a continuous-time control input $u_c:[0,\mathfrak{h}_{\mathfrak{S}_c}]\rightarrow\mathcal{U}$ that establishes \eqref{ControlProblem} if the MILP $\mathcal{Z}_{\mathfrak{S}_d}^{\mathrm{aug}}(u_d,X,U,z)$ is feasible. Moreover, for a sequence $u_d$ obtained from a feasible solution of $\mathcal{Z}_{\mathfrak{S}_d}^{\mathrm{aug}}(u_d,X,U,z)$, such continuous-time control input is given by \eqref{SpectralRepU}, where the matrix $U^i = [U_1^i \ U_2^i \ \ldots \ U_N^i]$ is defined according to \eqref{Partition} for all $i \in \llbracket 0,\mathfrak{h}_{\mathfrak{S}_d}\rrbracket$. 
\end{theorem}
\begin{proof}
	Suppose that the MILP $\mathcal{Z}_{\mathfrak{S}_d}^{\mathrm{aug}}(u_d,X,U,z)$ is feasible and let $u_d$ be the discrete-time control input obtained from a feasible solution. Bearing in mind that $z$ contains all binary variables recursively defined according to \eqref{PredicateConstraints}--\eqref{IntegratedConstraint}, we conclude that $\mathcal{Z}_{\mathfrak{S}_d}(u_d,z)$ is satisfied. This, by Proposition~\ref{PropMILPEncoding}, implies that \eqref{DiscreteControlProblem} holds. 
	
	After partitioning $X$ and $U$ according to \eqref{Partition}, we observe from $\mathcal{Z}_{\mathfrak{S}_d}^{\mathrm{aug}}(u_d,X,U,z)$ that \eqref{AlgorithmTechnical1} and \eqref{InputGuarantee} hold for all $i \in \llbracket 0,\mathfrak{h}_{\mathfrak{S}_d}-1\rrbracket$. We then construct $u_c\in\mathbb{I}_N^\tau(u_d)$ as in \eqref{SpectralRepU} and conclude from Theorem~\ref{ThControlConstraints} that it satisfies \eqref{DefInputConstraint}. Moreover, we note from Theorem~\ref{ThTrajectorySetCharacterization}
	that \eqref{SpectralRepX} holds. We will therefore show that \eqref{SpectralRepX} satisfies \eqref{ControlProblem}. 
	
	We accomplish this by utilizing structural induction to show that the satisfaction of each subformula in $\mathfrak{S}_c$ is implied by the satisfaction of its discrete-time counterpart in $\mathfrak{S}_d$. Specifically, for each subformula $\mathscr{S}_c$ in $\mathfrak{S}_c$ and its discrete-time counterpart $\mathscr{S}_d$ in $\mathfrak{S}_d$, we show that for all $i\in\llbracket 0, \mathfrak{h}_{\mathfrak{S}_d} \rrbracket$, the following implication holds:
	\begin{align}\label{SubformulaSatisfaction}
		&\Big(\restr{x_d(\cdot;x_0,u_d)}{\llbracket0,\mathfrak{h}_{\mathfrak{S}_d}\rrbracket},i\Big) \models \mathscr{S}_d\\ 
		&\Longrightarrow\forall t\in[0,\tau]: \Big(\restr{x_c(\cdot;x_0,u_c)}{[0,\mathfrak{h}_{\mathfrak{S}_c}]},i\tau+t\Big) \models \mathscr{S}_c.\nonumber
	\end{align}
	We note that in \eqref{SubformulaSatisfaction}, we adopt the standing convention that all time instants $i\tau+t$ are considered only for $t\in[0,\tau]$ such that $i\tau+t\leq \mathfrak{h}_{\mathfrak{S}_c}$, \textit{i.e.,} within the domain of $\restr{x_c(\cdot;x_0,u_c)}{[0,\mathfrak{h}_{\mathfrak{S}_c}]}$. This restriction is implicitly enforced throughout the induction. 
	
	We begin with the base case. Considering a predicate $\mu$ in $\mathfrak{S}_c$, let $\mathscr{S}_c = \mu$ and conclude from \eqref{SpecificationDiscretizationProcedure} that $\mathscr{S}_d = \mu$. If $\big(\restr{x_d(\cdot;x_0,u_d)}{\llbracket0,\mathfrak{h}_{\mathfrak{S}_d}\rrbracket},i\big) \models \mu$ for some $i\in\llbracket0,\mathfrak{h}_{\mathfrak{S}_d}\rrbracket$, it follows from \eqref{IntegratedConstraint} that \eqref{GuaranteeCondition} holds. This, by Theorem~\ref{ThGuaranteeCondition}, implies that $\big(\restr{x_c(\cdot;x_0,u_c)}{[0,\mathfrak{h}_{\mathfrak{S}_c}]},i\tau+t\big) \models \mu$ for all $t\in[0,\tau]$. We thus observe that \eqref{SubformulaSatisfaction} holds for all predicates $\mu$ in $\mathfrak{S}_c$. 
	
	 We proceed with the inductive hypothesis. Considering subformulas $\mathscr{S}_c^1$ and $\mathscr{S}_c^2$ in $\mathfrak{S}_c$, we let $\mathscr{S}_d^1$ and $\mathscr{S}_d^2$ be the corresponding discrete-time subformulas in $\mathfrak{S}_d$. We suppose that for all $i\in\llbracket0,\mathfrak{h}_{\mathfrak{S}_d}\rrbracket$, we have that 
	 \begin{align}
	 	&\Big(\restr{x_d(\cdot;x_0,u_d)}{\llbracket0,\mathfrak{h}_{\mathfrak{S}_d}\rrbracket},i\Big) \models \mathscr{S}_d^{l} \label{Assumption}\\
	 	&\Longrightarrow\forall t\in[0,\tau]: \Big(\restr{x_c(\cdot;x_0,u_c)}{[0,\mathfrak{h}_{\mathfrak{S}_c}]},i\tau+t\Big) \models \mathscr{S}_c^{l},\nonumber
	 \end{align}
	 where $l = 1,2$. 
	 
	 We now perform the inductive step. First, we let $\mathscr{S}_c$ be a Boolean combination (\textit{i.e.,} negation, conjunction, or disjunction) of $\mathscr{S}_c^1$ and $\mathscr{S}_c^2$. We observe from \eqref{SpecificationDiscretizationProcedure} that $\mathscr{S}_d$ is obtained by applying the same Boolean operator to $\mathscr{S}_d^1$ and $\mathscr{S}_d^2$. It is then straightforward to utilize the semantics \eqref{Semantics}, together with \eqref{Assumption}, to conclude \eqref{SubformulaSatisfaction}. 
	 
	 Next, we let $\mathscr{S}_c$ be a subformula of $\mathfrak{S}_c$ obtained by applying temporal operators (\textit{i.e.,} always, eventually, and until) to $\mathscr{S}_c^l$. We then note that $\mathscr{S}_d$ is obtained by applying the same temporal operators to $\mathscr{S}_d^l$. 
	 
	 For the always operator, let $\mathscr{S}_c= \operatorname{G}_{[t_1,t_2]}\mathscr{S}_c^1$ and recall from \eqref{SpecificationDiscretizationProcedure} that $\mathscr{S}_d = \operatorname{G}_{\llbracket\frac{t_1}{\tau},\frac{t_2}{\tau}\rrbracket}\mathscr{S}_d^1$. If $\mathscr{S}_d$ is valid at some instant $i\in\llbracket0,\mathfrak{h}_{\mathfrak{S}_d}\rrbracket$, it directly follows from \eqref{Semantics} that 
	 \begin{equation*}
	 	\forall i'\in\llbracket i+\frac{t_1}{\tau},i+\frac{t_2}{\tau}\rrbracket\!:\!  \Big(\restr{x_d(\cdot;x_0,u_d)}{\llbracket0,\mathfrak{h}_{\mathfrak{S}_d}\rrbracket},i' \Big) \!\models\! \mathscr{S}_d^1.
	 \end{equation*}
	 This, together with \eqref{Assumption}, implies that
	 \begin{align*}
	 \forall t\in[0,\tau]\!:\!  \Big(\restr{x_c(\cdot;x_0,u_c)}{[0,\mathfrak{h}_{\mathfrak{S}_c}]},i'\tau+t\Big)\! \models\! \mathscr{S}_c^1,
	 \end{align*}
	 for all $ i'\in\llbracket i+\frac{t_1}{\tau},i+\frac{t_2}{\tau}\rrbracket$, which indicates that $\mathscr{S}_c$ is valid at $i\tau+t$, for all $t\in[0,\tau]$. 
	  
	 Lastly, we consider the until operator and let $\mathscr{S}_c = \mathscr{S}_{c}^1\operatorname{U}_{[t_1,t_2]}\mathscr{S}_{c}^2$ and note that $\mathscr{S}_d= \mathscr{S}_{d}^1\operatorname{U}_{\llbracket\frac{t_1}{\tau},\frac{t_2}{\tau}\rrbracket}\mathscr{S}_{d}^2$. If 
	 $\mathscr{S}_d$ is valid at some instant $i\in\llbracket0,\mathfrak{h}_{\mathfrak{S}_d}\rrbracket$, 
	  we conclude from \eqref{Semantics} that there exists $i'\in\llbracket i+\frac{t_1}{\tau},i+\frac{t_2}{\tau}\rrbracket$ such that 
	 \begin{equation*}
	 	\forall i''\in \llbracket i,i'\rrbracket: \quad \Big(\restr{x_d(\cdot;x_0,u_d)}{\llbracket0,\mathfrak{h}_{\mathfrak{S}_d}\rrbracket},i''\Big) \models \mathscr{S}_{d}^1, 
	 \end{equation*}
	 and that 
	 \begin{equation*}
	 	\Big(\restr{x_d(\cdot;x_0,u_d)}{\llbracket0,\mathfrak{h}_{\mathfrak{S}_d}\rrbracket},i' \Big) \models \mathscr{S}_{d}^2.
	 \end{equation*}
	 This, together with \eqref{Assumption}, implies that 
	 \begin{equation*}
	 	\forall t\in [0,\tau]: \quad \Big(\restr{x_c(\cdot;x_0,u_c)}{[0,\mathfrak{h}_{\mathfrak{S}_c}]},i''\tau+t\Big) \models \mathscr{S}_{c}^1, 
	 \end{equation*}
	 for all $i''\in \llbracket i,i'\rrbracket$, and that 
	 \begin{equation*}
	 	\forall t\in[0,\tau]: \quad \Big(\restr{x_c(\cdot;x_0,u_c)}{[0,\mathfrak{h}_{\mathfrak{S}_c}]},i'\tau+t \Big) \models \mathscr{S}_{c}^2,
	 \end{equation*}
	 indicating that $\mathscr{S}_c$ is valid at $i\tau+t$, for all $t\in[0,\tau]$. It is worth mentioning that the subformulas formed by the (eventually) operator $\operatorname{F}$ do not need to be addressed separately, as they can be expressed in terms of the (until) operator $\operatorname{U}$. We therefore conclude that \eqref{SubformulaSatisfaction} holds for any subformula $\mathscr{S}_c$ obtained by applying temporal operators (\textit{i.e.,} always, eventually, and until) on $\mathscr{S}_c^l$.
	 
	 Having established \eqref{SubformulaSatisfaction} for any subformula $\mathscr{S}_c$ in $\mathfrak{S}_d$, we apply induction on the structure of $\mathfrak{S}_d$ to conclude that 
	 \begin{align*}
	 	&\restr{x_d(\cdot;x_0,u_d)}{\llbracket0,\mathfrak{h}_{\mathfrak{S}_d}\rrbracket} \models \mathfrak{S}_d\Longrightarrow\restr{x_c(\cdot;x_0,u_c)}{[0,\mathfrak{h}_{\mathfrak{S}_c}]} \models \mathfrak{S}_c.
	 \end{align*}
	This, together with \eqref{DiscreteControlProblem} and \eqref{Semantics}, implies that \eqref{ControlProblem} holds.
\end{proof}

Theorem~\ref{ThMain} characterizes the continuous-time STL control problem as the (augmented) MILP $\mathcal{Z}_{\mathfrak{S}_d}^{\mathrm{aug}}(u_d,X,U,z)$. A feasible solution to this MILP gives the control sequence $u_d$, based on which one may obtain the appropriate continuous-time control input $u_c$ as in \eqref{SpectralRepU}. Based on Theorem~\ref{ThMain}, we may summarize the hierarchical continuous-time STL control scheme in Algorithm~\ref{Algorithm}.
\begin{algorithm}
	\caption{Hierarchical continuous-time STL control}\label{Algorithm}
	\begin{algorithmic}[1]
		\Require $A_c$, $B_c$, $\mathfrak{S}_c$, $N$
		\algrule
		\vskip.2\baselineskip
		\Statex \hspace*{-4.6ex} \textbf{1. Discretization}
		\vskip.5\baselineskip
		\State Select the sampling time $\tau>0$ such that every time bound specified in $\mathfrak{S}_c$ is an integer multiple of $\tau$;
		\State Apply the recursive procedure \eqref{SpecificationDiscretizationProcedure} to construct $\mathfrak{S}_d$; 
		\State Define the matrices $Q_N^\tau(A_c,B_c)$, $R_N^\tau(A_c,B_c)$, $T_1$ , and $T_2$ as in \eqref{QRSV};
		\State Solve \eqref{DesignResult} and obtain the matrices $A_d$ and $B_d$ according to \eqref{IntermediateResult};
		\algrule
		\setcounter{ALG@line}{0}
		\algrenewcommand\alglinenumber[1]{\footnotesize2.#1:}
		\vskip.2\baselineskip
		\Statex \hspace*{-4.6ex} \textbf{2. Planning}
		\vskip.5\baselineskip
		\State For each $i\in\llbracket0,\mathfrak{h}_{\mathfrak{S}_d}\rrbracket$, define the equation \eqref{AlgorithmTechnical1} and the inclusion \eqref{InputGuarantee};
		\State For each predicate \eqref{Predicate} in $\mathfrak{S}_d$, define the inequality \eqref{IntegratedConstraint};
		\State Introduce the variables $X$ and $U$ as in \eqref{Partition};
		\State Obtain the set of constraints $\mathcal{Z}_{\mathfrak{S}_d}(u_d,z)$ according to \eqref{PredicateConstraints}--\eqref{SatisfactionConstraint};
		\State Collect the constraints from \eqref{AlgorithmTechnical1}, \eqref{InputGuarantee}, and $\mathcal{Z}_{\mathfrak{S}_d}(u_d,z)$ into the augmented MILP $\mathcal{Z}_{\mathfrak{S}_d}^{\text{aug}}(u_d,X,U,z)$;
		\State Solve $\mathcal{Z}_{\mathfrak{S}_d}^{\text{aug}}(u_d,X,U,z)$ and obtain $u_d$;
		\algrule
		\setcounter{ALG@line}{0}
		\algrenewcommand\alglinenumber[1]{\footnotesize3.#1:}
		\vskip.2\baselineskip
		\Statex \hspace*{-4.6ex} \textbf{3. Continuous-time Execution}
		\vskip.5\baselineskip
		\State Construct the control input $u_c:[0,\mathfrak{h}_{\mathfrak{S}_c}]\rightarrow\mathcal{U}$ as in \eqref{SpectralRepU};
		\algrule
		\Ensure The continuous-time control input $u_c$
	\end{algorithmic}
\end{algorithm} 

\section{Numerical Example}\label{Sec_Example}
We consider the STL control of a robotic agent $\bm{\Sigma}_c$, whose motion is described by 
\begin{equation*}
	\begin{aligned}
		\dot{x}_{\mathrm{cont}} (t) &= v_{\mathrm{cont}}(t);\\
		\dot{v}_{\mathrm{cont}} (t) &= F_{\mathrm{cont}}(t),
	\end{aligned}
\end{equation*}
where the variables $x_{\mathrm{cont}}$, $v_{\mathrm{cont}}$, and $F_{\mathrm{cont}}$ respectively denote displacement, velocity, and acting force in the longitudinal direction. Due to physical limitations, the acting force must not exceed $200 \unit{\meter/\second^2}$, \textit{i.e.,} $|F_{\mathrm{cont}}| \leq 200$. We suppose that the robot $\bm{\Sigma}_c$ is initially at rest and located at $1$, \textit{i.e.,} $x_{\mathrm{cont}}(0) = 1$ and $v_{\mathrm{cont}}(0) = 0$. The objective is to apply the acting force $F_{\mathrm{cont}}$ such that the robot $\bm{\Sigma}_c$ evolves according to the STL formula
\begin{equation}\label{Sc}
	\begin{aligned}
		\mathfrak{S}_c = & \ \mathrm{F}_{[0.2,0.8]}\big(x_{\mathrm{cont}}(t) \leq -1.5 \big) \\
		&\wedge \mathrm{F}_{[1,1.4]}\big(x_{\mathrm{cont}}(t) \geq 1.5 \big)\\
		& \wedge \mathrm{F}_{[1.6,2]}\big(x_{\mathrm{cont}}(t) \leq -1.5 \big)\\
		& \wedge \mathrm{G}_{[0,2]} \big(v_{\mathrm{cont}}(t) \leq 18 \wedge v_{\mathrm{cont}}(t) \geq -18\big).
	\end{aligned}
\end{equation}
We note that the STL formula $\mathfrak{S}_c$ requires the robot to 1) \emph{eventually} reach the location less than or equal to $-1.5 \unit{\meter}$ between $0.2$ and $0.8$ seconds, 2) \emph{eventually} reach the location greater than or equal to $1.5 \unit{\meter}$ between $1$ and $1.4$ seconds, 3) \emph{eventually} reach the location less than or equal to $-1.5 \unit{\meter}$ between $1.6$ and $2$ seconds, and 4) \emph{always} maintain its velocity between $-18 \unit{\meter/\second}$ and $18 \unit{\meter/\second}$. It is also worth mentioning that $\mathfrak{h}_{\mathfrak{S}_c} = 2$. 

We utilize the proposed hierarchical continuous-time STL control scheme to obtain the acting force that enforces $\bm{\Sigma}_c$ to evolve in accordance with the specification $\mathfrak{S}_c$. 

Taking $x_c = (x_{\mathrm{cont}},v_{\mathrm{cont}})^\top$ and $u_c = F_{\mathrm{cont}}$, we formulate $\bm{\Sigma}_c$ as in \eqref{ContinuousSystem} with 
\begin{equation*}
	\begin{aligned}
		A_c &= \begin{bmatrix}
			0&1\\
			0&0
		\end{bmatrix}, &B_c&=\begin{bmatrix}
			0\\
			1
		\end{bmatrix}.
	\end{aligned}
\end{equation*}
Bearing in mind that $|F_{\mathrm{cont}}| \leq 200$, we formulate the control input constraints by restricting $u_c$ to take values in the convex polyhedral set 
\begin{equation*}
	\mathcal{U} = \left\{u\in\mathbb{R}^m \left \vert \begin{bmatrix}
		1\\
		-1
	\end{bmatrix} \right. u \leq \begin{bmatrix}
	200\\
	200
	\end{bmatrix}\right\}.
\end{equation*}
We now proceed in accordance with Algorithm~\ref{Algorithm} to hierarchically construct the control input $u_c$ that establishes \eqref{ControlProblem}. 

First, we choose the sampling time $\tau$ such that every time bound specified in $\mathfrak{S}_c$ appears as an integer multiple of $\tau$. We observe from \eqref{Sc} that one such choice is $\tau = 0.2$.

Second, we apply the recursive procedure \eqref{SpecificationDiscretizationProcedure} to obtain the discrete-time STL formula
\begin{equation*}
	\begin{aligned}
		\mathfrak{S}_d =& \ \mathrm{F}_{\llbracket 1,4\rrbracket}\big(x_{\mathrm{disc}}(i) \leq -1.5 \big) \\
		&\wedge \mathrm{F}_{\llbracket 5,7\rrbracket}\big(x_{\mathrm{disc}}(i) \geq 1.5 \big)\\
		& \wedge \mathrm{F}_{\llbracket 8,10\rrbracket}\big(x_{\mathrm{disc}}(i) \leq -1.5 \big)\\
		& \wedge \mathrm{G}_{\llbracket 0,10\rrbracket} \big(v_{\mathrm{disc}}(i) \leq 18 \wedge v_{\mathrm{disc}}(i) \geq -18\big),
	\end{aligned}
\end{equation*}
where $x_{\mathrm{disc}}$ and $v_{\mathrm{disc}}$ are respectively discrete-time counterparts to $x_{\mathrm{cont}}$ and $v_{\mathrm{cont}}$. It is worth noting that $\mathfrak{h}_{\mathfrak{S}_d} = 10$.

Third, we select the polynomial degree $N = 5$ and define the matrices $Q_N^\tau(A_c,B_c)$, $R_N^\tau(A_c,B_c)$, $T_1$ , and $T_2$ as in \eqref{QRSV}. Accordingly, we utilize Theorem~\ref{Th_Design} to construct the discrete-time system $\bm{\Sigma}_d$ as in \eqref{DiscreteSystem} with 
\begin{equation*}
	\begin{aligned}
		A_d&= \begin{bmatrix}
			0.6990&0.1398\\
			0.0000&0.6990
		\end{bmatrix}, &B_d&=\begin{bmatrix}
			0.0000\\0.1398
		\end{bmatrix}.
	\end{aligned}
\end{equation*} 
Hence, $\bm{\Sigma}_c$ is an $N$-th order interpolator of $\bm{\Sigma}_d$ with respect to the sampling time $\tau$. We also note that the matrices $A_d$ and $B_d$ are not necessarily unique. 

Fourth, we recursively construct the set of constraints $\mathcal{Z}_{\mathfrak{S}_d}^{\mathrm{aug}}(u_d,X,U,z)$ according to \eqref{AlgorithmTechnical1}, \eqref{InputGuarantee}, and \eqref{PredicateConstraints}--\eqref{IntegratedConstraint}. We then employ the optimization toolbox of MATLAB 2024b to solve the MILP $\mathcal{Z}_{\mathfrak{S}_d}^{\mathrm{aug}}(u_d,X,U,z)$ and, accordingly, obtain the control sequence $u_d:\llbracket0,\mathfrak{h}_{\mathfrak{S}_d}\rrbracket\rightarrow\mathcal{U}$ that establishes \eqref{DiscreteControlProblem}. 

At last, we construct the control input $u_c:[0,\mathfrak{h}_{\mathfrak{S}_c}]\rightarrow\mathcal{U}$ as in \eqref{SpectralRepU} and conclude from Theorem~\ref{ThMain} that \eqref{ControlProblem} holds. 

\begin{figure}
	\centering
	\begin{subfigure}[b]{0.5\textwidth}
		\includegraphics[width=1\textwidth]{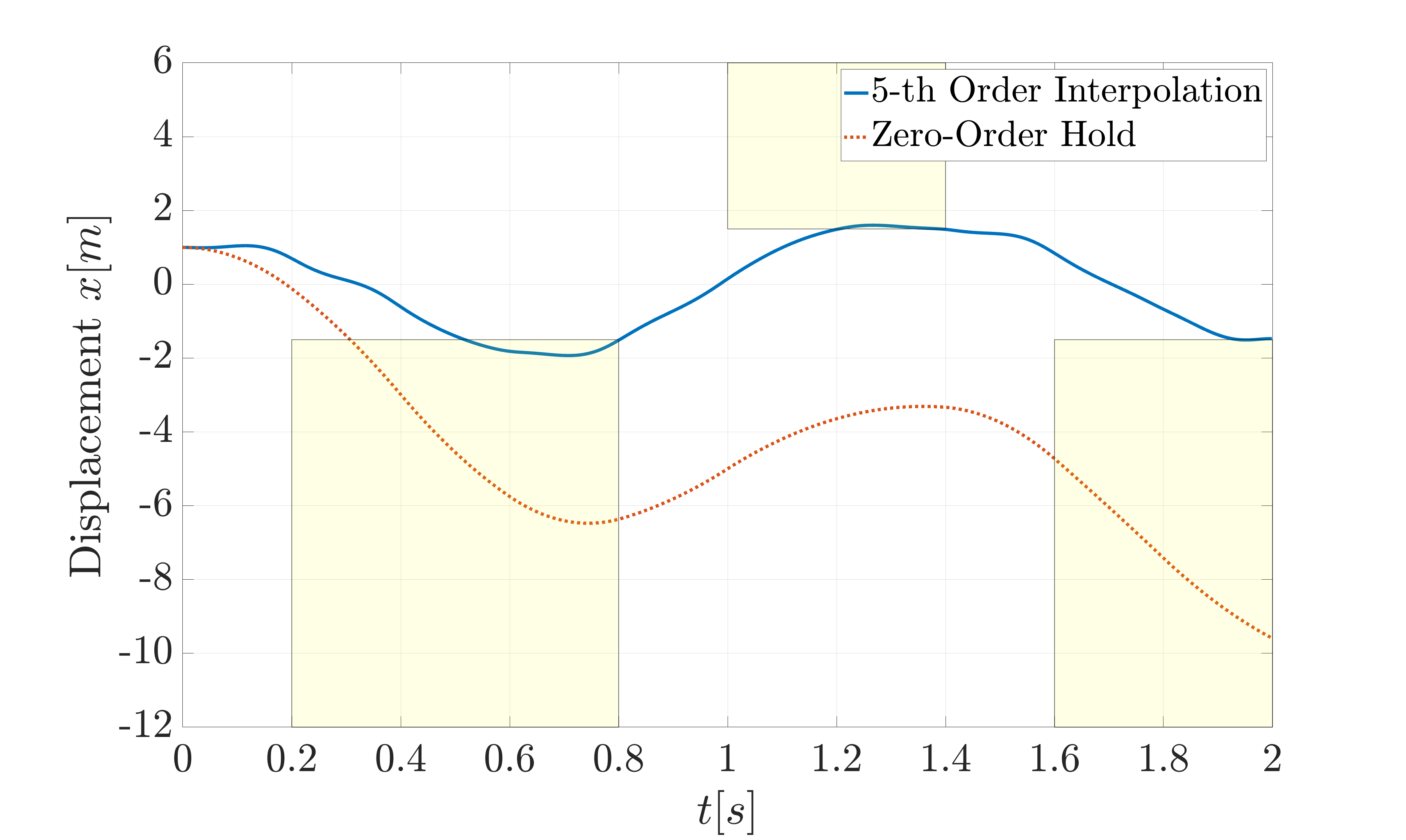}
		\caption{Displacement $x$.}
		\label{Displacement}
	\end{subfigure}
	\vfill
	\begin{subfigure}[b]{0.5\textwidth}
		\includegraphics[width=1\textwidth]{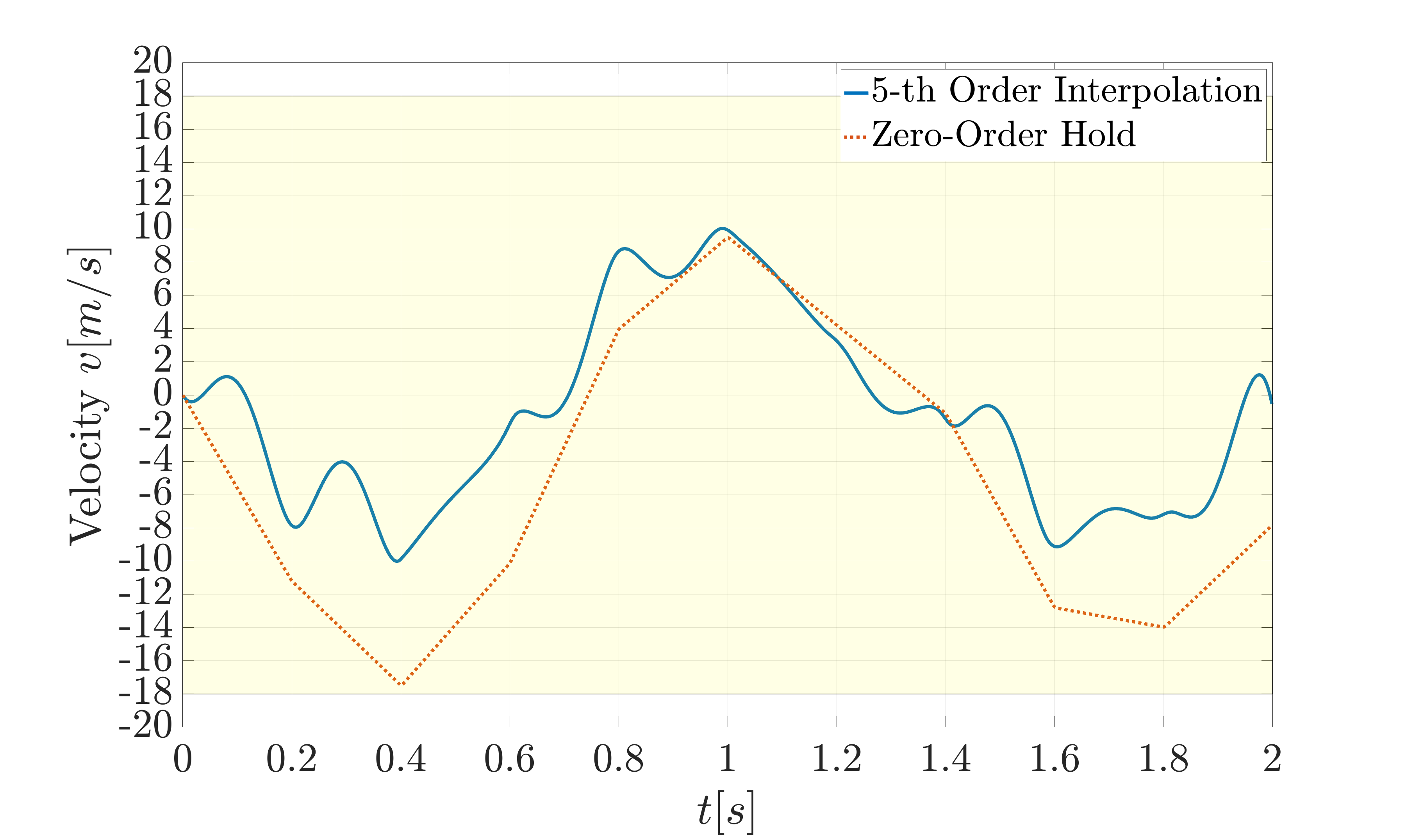}
		\caption{Velocity $v$.}
		\label{Velocity}
	\end{subfigure}
	\caption{Subject to the continuous-time control \eqref{SpectralRepU}, which is indicated by solid lines in Figure~\ref{Inputs}, the robot $\bm{\Sigma}_c$ experiences a displacement and a velocity that comply with the requirements specified by the STL formula $\mathfrak{S}_c$. By contrast, subject to a continuous-time control $u_c$ obtained by ZOH interpolation (indicated by dashed lines in Figure~\ref{Inputs}), the evolution of $\bm{\Sigma}_c$ does not adhere to the formula $\mathfrak{S}_c$.}
	\label{Outputs}
\end{figure}

\begin{figure}
	\centering
	\includegraphics[width=0.5\textwidth]{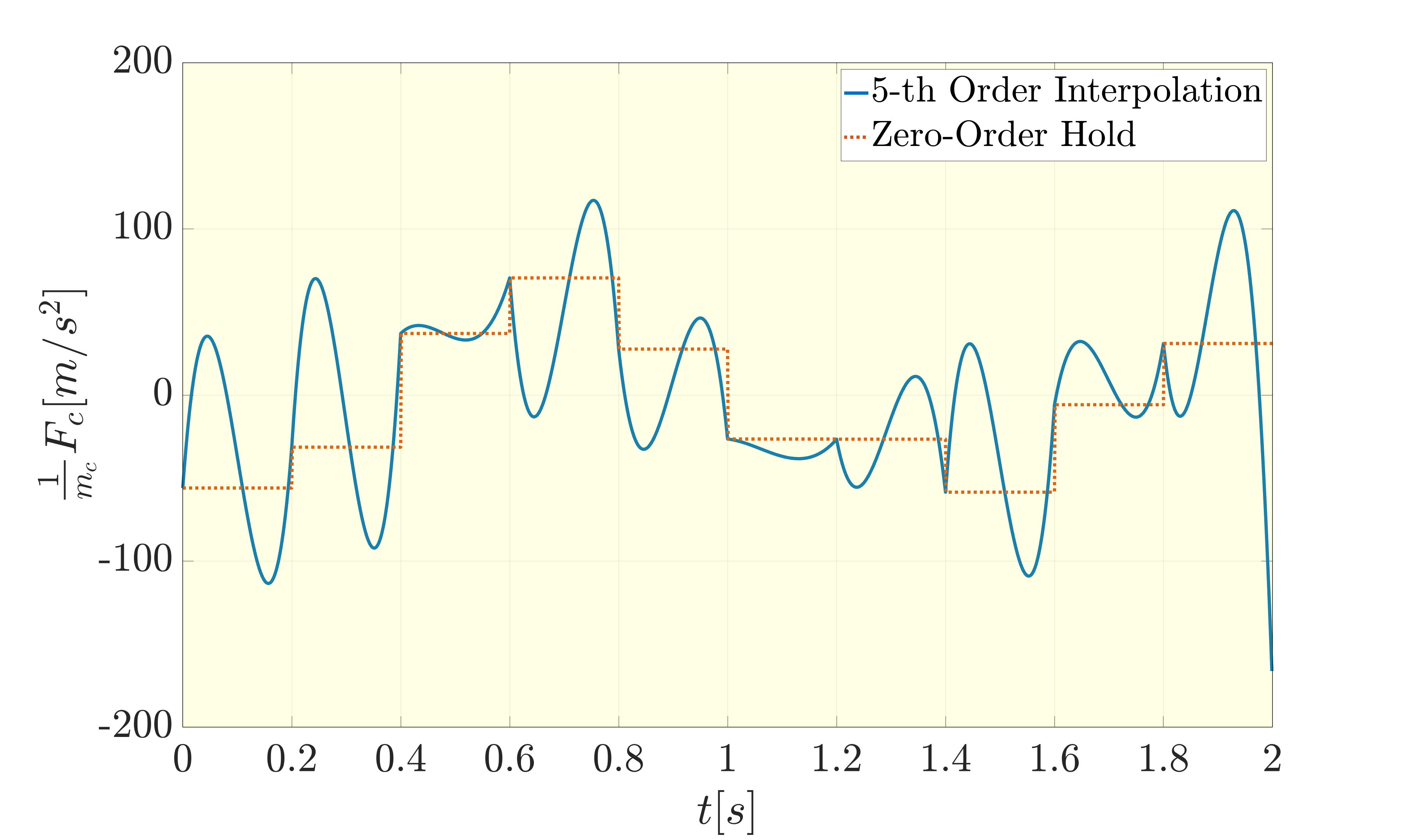}
	\caption{The continuous-time control \eqref{SpectralRepU}, which is indicated by solid lines, enforces the robot $\bm{\Sigma}_c$ to evolve in accordance with the STL formula $\mathfrak{S}_c$ (Figure~\ref{Velocity}). Importantly, such control respects the input constraints as its magnitude does not exceed $200 \unit{\meter/\second^2}$. By contrast, subject to the continuous-time control $u_c$ obtained by ZOH interpolation of the sequence $u_d$ (indicated by dashed lines), the evolution of $\bm{\Sigma}_c$ does not comply with formula $\mathfrak{S}_c$.}
	\label{Inputs}
\end{figure}

The displacement $x_{\mathrm{cont}}$ and the velocity $v_{\mathrm{cont}}$ of the robot $\bm{\Sigma}_c$ are depicted in Figure~\ref{Outputs} (indicated by solid lines), whereas the continuous-time control \eqref{SpectralRepU} is depicted in Figure~\ref{Inputs} (indicated by solid lines). In each figure, the regions of interest (specified by the STL formula $\mathfrak{S}_c$) are indicated by colored patches. It is clear from Figure~\ref{Outputs} that the continuous-time control \eqref{SpectralRepU} establishes \eqref{ControlProblem} while respecting the input constraints (as illustrated in Figure~\ref{Inputs}), \textit{i.e.,} the acting force enforces the robot to evolve in accordance with the specification $\mathfrak{S}_c$ while its magnitude does not exceed $200 \unit{\meter/\second^2}$. 

We show the efficacy of our hierarchical scheme by comparing the performance of the control input  \eqref{SpectralRepU} with that of an input obtained as a zero-order hold (ZOH) interpolation of the sequence $u_d$ (indicated by dashed lines in Figure~\ref{Inputs}). The displacement $x_{\mathrm{cont}}$ and the velocity $v_{\mathrm{cont}}$ of the robot $\bm{\Sigma}_c$, subject to this control, are also depicted (in dashed lines) in Figure~\ref{Outputs}. We observe that subject to such control, the robot does not fulfill the requirements specified by the formula $\mathfrak{S}_c$. Particularly, in this case, the displacement $x_{\mathrm{cont}}$ never reaches $1.5 \unit{m}$ between $1$ and $1.4$ seconds. 

\section{Conclusion}\label{Sec_Conclusion}
We proposed a hierarchical scheme that solves the continuous-time STL control problem by translating it into a discrete-time one in such a way that solutions to the original continuous-time problem can be characterized in terms of solutions to the discrete-time problem. We implemented this hierarchical schemes according to the following steps. First, we discretized the continuous-time STL formula into a discrete-time one by sampling the temporal behavior it specifies. Accordingly, we discretized the continuous-time system into a discrete-time model such that the input-state trajectories of the continuous-time system can be realized as piecewise polynomial interpolations of the input-state trajectories of the discrete-time model. We accomplished this by utilizing the notion of system interpolation, which enabled us to characterize the input-state behavior of the continuous-time system in terms of its discrete-time model. Second, we exploited this characterization to formulate mixed-integer linear constraints on the input-state trajectories of the discrete-time model and, accordingly, plan a trajectory such that its piecewise polynomial interpolation satisfies the continuous-time STL formula. As a matter of fact, we implemented the discretization and planning procedures in such a manner that any discrete-time trajectory that respects these constraints can be transformed into a continuous-time trajectory that is guaranteed to satisfy the original STL formula. Lastly, we proposed a step-by-step algorithm that constructs the control input that realizes the such continuous-time trajectory. 

For the future, we aim to extend the proposed hierarchical continuous-time STL control scheme to nonlinear systems.  
\appendix
\begin{lemma}\label{TechnicalLemma}
	The matrix $\mathfrak{B}$, defined in \eqref{B}, is invertible.
\end{lemma}
\begin{proof}
	We note that since \eqref{B} is a square matrix, its invertibility is equivalent to its non-singularity. Hence, in order to show that \eqref{B} is invertible, it suffices to show that it is non-singular. To this end, we suppose that $\mathfrak{B}$ is singular. This implies the existence of a non-zero vector $\nu \in \mathbb{R}^{N+1}$ such that $\mathfrak{B}\nu = 0$. We let $\nu_0,\nu_1,\ldots,\nu_N\in\mathbb{R}$ be such that $\nu = (\nu_0,\nu_1,\cdots,\nu_N)^\top$ and define the polynomial $q\in\restr{P_N^1}{[0,\tau]}$ such that 
	\begin{equation}\label{Technical1}
		\forall t\in[0,\tau]: \quad q(t) = \sum_{k=0}^{N}\nu_k\mathfrak{b}_{k,N}\left(\frac{t}{\tau}\right).
	\end{equation}
	This, as a consequence of fact that $\mathfrak{B}\nu = 0$, implies that 
	\begin{equation}\label{Technical2}
		\forall j\in\llbracket0,N\rrbracket: \quad q(t_{N,j}) = 0.
	\end{equation} 
	We then recall that $q\in\restr{P_N^1}{[0,\tau]}$ and conclude from \eqref{SpectralRepresentation} and \eqref{Technical2} that 
	\begin{equation*}
		\forall t\in[0,\tau]: \quad q(t) = 0.
	\end{equation*}
	It then follows from \eqref{Technical1} that
	\begin{equation*}
		\sum_{k=0}^{N}\nu_k\mathfrak{b}_{k,N}\left(\frac{t}{\tau}\right) = 0,
	\end{equation*}
	which, due to the linear independence of the Bernstein basis polynomials, implies that $\nu = 0$, contradicting the initial assumption that $\nu$ is non-zero. We thus conclude, by contradiction, that $\mathfrak{B}$ is non-singular. 
\end{proof}
\bibliographystyle{ieeetr}
\bibliography{Reference.bib}
\newpage
\begin{IEEEbiography}[{\includegraphics[width=1in,height=1.25in,clip,keepaspectratio]{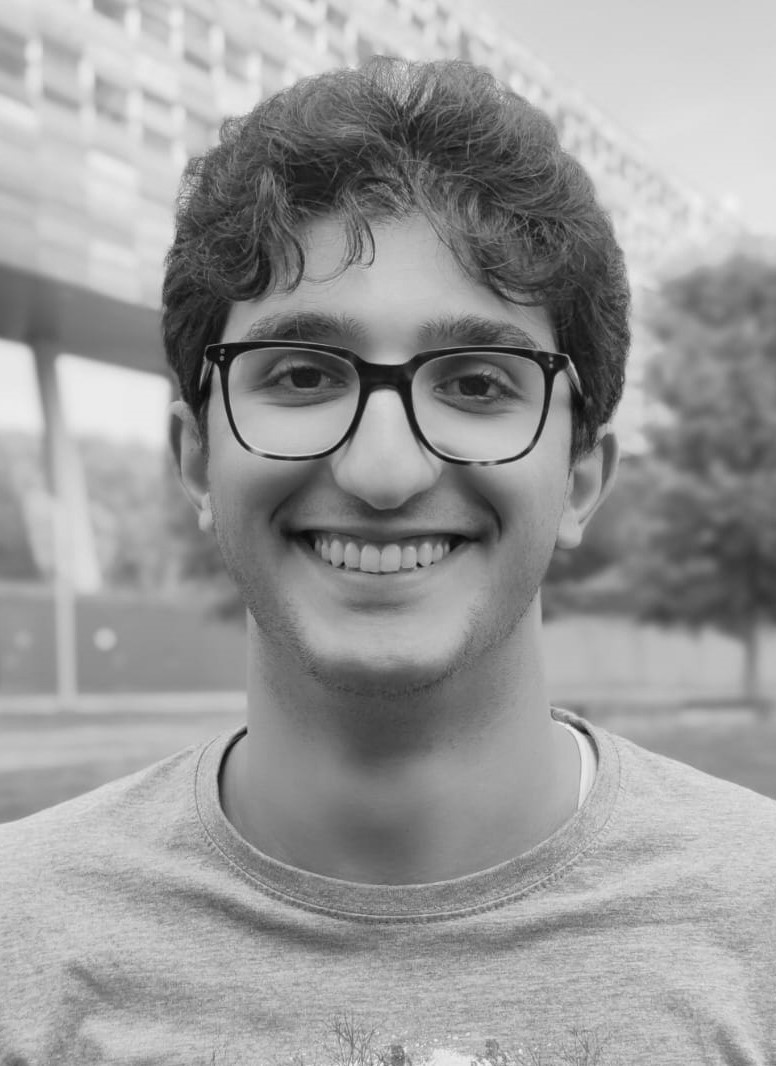}}]{Armin Pirastehzad} (Member, IEEE) earned his Ph.D. in Applied Mathematics, specializing in Systems and Control, from the University of Groningen, the Netherlands, in 2025. He previously obtained his B.Sc. and M.Sc. degrees in Electrical Engineering, also specializing in Systems and Control, from the University of Tehran, Iran, in 2017 and 2020, respectively. He is currently a Postdoctoral Researcher at the Bernoulli Institute for Mathematics, Computer Science, and Artificial Intelligence, University of Groningen.
\end{IEEEbiography}
\vspace*{-10cm}
\begin{IEEEbiography}[{\includegraphics[width=1in,height=1.25in,clip,keepaspectratio]{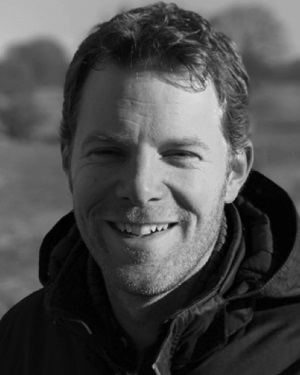}}]{Bart Besselink} (Senior Member, IEEE) received the M.Sc. (cum laude) degree in mechanical engineering in 2008 and the Ph.D. degree in 2012, both from Eindhoven University of Technology, Eindhoven, The Netherlands.
	
	Since 2016, he has been with the Bernoulli Institute for Mathematics, Computer Science and Artificial Intelligence, University of Groningen, Groningen, The Netherlands, where he is currently an associate professor. He was a short-term Visiting Researcher with the Tokyo Institute of Technology, Tokyo, Japan, in 2012. Between 2012 and 2016, he was a Postdoctoral Researcher with the ACCESS Linnaeus Centre and Department of Automatic Control, KTH Royal Institute of Technology, Stockholm, Sweden.
	
	His main research interests are on mathematical systems theory for large-scale interconnected systems, with emphasis on contract-based design and control, compositional analysis, model reduction, and applications in intelligent transportation systems and neuromorphic computing. He is a recipient (with Xiaodong Cheng and Jacquelien Scherpen) of the 2020 Automatica Paper Prize.
\end{IEEEbiography}

\end{document}